\documentclass[12pt, twoside]{article}
\usepackage{times}
\usepackage{graphics}
\usepackage{theorem}
\usepackage{graphicx}
\usepackage{amsmath}
\usepackage{latexsym}
\usepackage{amssymb}
\usepackage[flushmargin]{footmisc}
\renewcommand{\thefootnote}{\fnsymbol{footnote}}
\renewcommand{\baselinestretch}{1.1}
\makeatletter
\def\eqnarray{\stepcounter{equation}\let\@currentlabel=\theequation
\global\@eqnswtrue
\global\@eqcnt\z@\tabskip\@centering\let\\=\@eqncr
$$\halign to \displaywidth\bgroup\@eqnsel\hskip\@centering
  $\displaystyle\tabskip\z@{##}$&\global\@eqcnt\@ne 
  \hfil$\;{##}\;$\hfil
  &\global\@eqcnt\tw@ $\displaystyle\tabskip\z@{##}$\hfil 
   \tabskip\@centering&\llap{##}\tabskip\z@\cr}
\makeatother
\theorembodyfont{\itshape}

\newtheorem{thm}{Theorem}[section]
\newtheorem{lem}{Lemma}[section]
\newtheorem{prop}{Proposition}[section]
\newtheorem{coro}{Corollary}[section]

\theorembodyfont{\rmfamily}
\newenvironment{proof}{\noindent{\it Proof.~~}}{\qed}

\newtheorem{rem}{Remark}[section]{\it}{\rm}
\newtheorem{defn}{Definition}[section]{\bf}{\rm}
\newtheorem{assumption}{Assumption}[section]{\bf}{\it}

\def\vc#1{\mbox{\boldmath $#1$}}
\def\svc#1{\mbox{\boldmath $\scriptstyle #1$}}

\newcommand{\qed}{\hspace*{\fill}$\Box$}
\newcommand{\dm}{\displaystyle}
\newcommand{\diag}{\mathrm{diag}}
\newcommand{\conj}{\mathrm{conj}}
\newcommand{\adj}{\mathrm{adj}}
\newcommand{\Mod}{\mathrm{mod}}

\newcommand{\rmd}{\mathrm{d}}
\newcommand{\E}{\mathrm{E}}
\newcommand{\bbH}{\mathbb{H}}
\newcommand{\bbL}{\mathbb{L}}
\newcommand{\J}{\mathbb{J}}
\newcommand{\K}{\mathbb{K}}

\newcommand{\M}{\mathbb{M}}
\newcommand{\N}{\mathbb{N}}
\newcommand{\R}{\mathbb{R}}

\newcommand{\calP}{\mathcal{P}}
\newcommand{\calA}{\mathcal{A}}

\newcommand{\Z}{\mathbb{Z}}
\newcommand{\C}{\mathbb{C}}
\newcommand{\D}{\mathbb{D}}

\makeatother

\begin{document}\thispagestyle{plain} 

\hfill

{\Large{\bf
\begin{center}
Light-tailed asymptotics of stationary tail probability vectors of Markov chains of M/G/1 type%
\footnote[1]{This is a revised version of the paper published in Stochastic Models vol.~26, no.~4, pp.~505--548, 2010. \\
In the revised version, some errors are corrected.}  
\end{center}
}
}

\begin{center}
{
Tatsuaki Kimura, Kentaro Daikoku, Hiroyuki Masuyama%
\footnote[2]{Address correspondence to Hiroyuki Masuyama, E-mail:
masuyama@sys.i.kyoto-u.ac.jp},
and Yutaka Takahashi
}

\medskip

{\small
Department of Systems
Science, Graduate School of Informatics, Kyoto University\\
Kyoto 606-8501, Japan
}

\bigskip
\medskip

{\small 
\textbf{Abstract}

\medskip

\begin{tabular}{p{0.85\textwidth}}
This paper studies the light-tailed asymptotics of the stationary tail
probability vectors of a Markov chain of M/G/1 type.  Almost all
related studies have focused on the {\it typical} case, where the
transition block matrices in the non-boundary levels have a dominant
impact on the decay rate of the stationary tail probability vectors
and their decay is aperiodic. In this paper, we study not only the
typical case but also {\it atypical} cases such that the stationary
tail probability vectors decay periodically and/or their decay rate is
determined by the tail distribution of jump sizes from the boundary
level. We derive light-tailed asymptotic formulae for the stationary
tail probability vectors by locating the {\it dominant poles} of the
generating function of the sequence of those vectors. Further we
discuss the positivity of the dominant terms of the obtained
asymptotic formulae.
\end{tabular}
}
\end{center}

\begin{center}
\begin{tabular}{p{0.90\textwidth}}
{\small

\medskip

}
\end{tabular}

\end{center}

\section{Introduction}\label{introduction}
\renewcommand{\thefootnote}{$\ddag$\arabic{footnote}}

This paper considers a Markov chain $\{(X_n,S_n); \; n=0,1,\dots\}$ of
M/G/1 type \cite{Neut89}, where $X_n \in \{0,1,\dots\}$ and
\[
\begin{array}{ll}
S_n \in \M_0 \triangleq \{1,2,\dots,M_0\}, & \mbox{if}~X_n = 0,
\\
S_n \in \M \triangleq \{1,2,\dots,M\}, & \mbox{otherwise}.
\end{array}
\]
The sets of states $\{(0,j); j \in \M_0\}$ and $\{(k,j); j \in \M\}$
($k=1,2,\dots$) are called level 0 and level $k$, respectively.
Arranging the states in lexicographical order, the transition
probability matrix $\vc{T}$ of $\{(X_n,S_n); \; n=0,1,\dots\}$ is
given by
\begin{equation}
\vc{T} = 
\bordermatrix{
               & {\rm lev.\ 0}&   1  &  2       &  3       & \cdots &       
\cr
{\rm lev.\ 0}\ & \vc{B}(0) & \vc{B}(1) & \vc{B}(2) & \vc{B}(3) & \cdots
\cr
\phantom{\rm lev.\ }1 & \vc{C}(0) & \vc{A}(1) & \vc{A}(2) & \vc{A}(3) & \cdots
\cr
\phantom{\rm lev.\ }2 & \vc{O}   & \vc{A}(0) & \vc{A}(1) & \vc{A}(2) & \cdots
\cr
\phantom{\rm lev.\ }3 & \vc{O}   & \vc{O}   & \vc{A}(0) & \vc{A}(1) & \cdots
\cr
\phantom{\rm lev.\ }\vdots    & \vdots   & \vdots  &  \vdots   & \vdots   & \ddots
},
\label{eqn-matrix-T}
\end{equation}
where $\vc{A}(k)$ ($k=0,1,\dots$) is an $M \times M$ matrix,
$\vc{B}(0)$ is an $M_0 \times M_0$ matrix, $\vc{B}(k)$ ($k=1,2,\dots$)
is an $M_0 \times M$ matrix, and $\vc{C}(0)$ is an $M \times M_0$
matrix.

We define $\vc{A}$ and $\vc{B}$ as
\[
\vc{A}
=
\sum_{k=0}^{\infty}
\vc{A}(k),
\qquad
\vc{B}
=
\sum_{k=1}^{\infty}
\vc{B}(k),
\]
respectively. We assume that $\vc{A}$ is a stochastic matrix and
$\vc{B}(0)\vc{e}+\vc{B}\vc{e}=\vc{e}$, where $\vc{e}$ denotes a column
vector of ones with an appropriate dimension. Let $\vc{\pi}$ denote a
$1 \times M$ vector such that $\vc{\pi}\vc{A} = \vc{\pi}$ and
$\vc{\pi}\vc{e} = 1$. Note here that if $\vc{A}$ is irreducible,
$\vc{\pi}$ is uniquely determined.

Throughout this paper, we make the following assumption.
\begin{assumption}\label{assu-1}
\hfill
\begin{enumerate}
\item $\vc{T}$ is irreducible.

\item
$\vc{A}$ is irreducible.

\item $\rho \triangleq \vc{\pi} \vc{\beta}_A < 1$, where $\vc{\beta}_A
= \sum_{k=1}^{\infty} k \vc{A}(k) \vc{e}$.
\item
$\vc{\beta}_B \triangleq \sum_{k=1}^{\infty}
k \vc{B}(k) \vc{e}$ is a finite vector.
\end{enumerate}
\end{assumption}

Let $\vc{x}$ denote the stationary probability vector of $\vc{T}$,
i.e., $\vc{x}\vc{T} = \vc{x}$ and $\vc{x}\vc{e}=1$. It is known that
under Assumption \ref{assu-1}, $\{(X_n,S_n); \; n=0,1,\dots\}$ is
irreducible and positive recurrent (see Proposition~3.1 in Chapter~XI
of \cite{Asmu03}). Therefore if Assumption \ref{assu-1} holds, then
$\vc{x} > \vc{0}$, which is uniquely determined.  Let $\vc{x}(k)$
($k=0,1,\dots$) denote a subvector of $\vc{x}$ corresponding to level
$k$. We then have $\vc{x} = (\vc{x}(0),\vc{x}(1),\vc{x}(2),\dots)$ and
\begin{equation}
\vc{x}(k) 
= \vc{x}(0)\vc{B}(k) + \sum_{l=1}^{k+1}\vc{x}(l)\vc{A}(k+1-l),
\qquad k=1,2,\dots.
\label{eqn-x(k)}
\end{equation}
Further let $\overline{\vc{x}}(k)=\sum_{l=k+1}^{\infty}\vc{x}(l)$
($k=0,1,\dots$), which is a positive vector. We call
$\overline{\vc{x}}(k)$'s stationary tail probability vectors of
$\vc{T}$ hereafter.

In this paper, we study the light-tailed asymptotics of
$\{\overline{\vc{x}}(k); k=0,1,\dots\}$. The following is the
definition of a light-tailed sequence of nonnegative matrices
(including vectors).
\begin{defn}\label{def-light-tailed}
Let $\vc{Y}(k)$'s ($k=0,1,\dots$) denote nonnegative matrices with the
same dimension, and let $\vc{Y}^{\ast}(z)$ denote the generating
function of $\{\vc{Y}(k)\}$ defined by the power series
$\sum_{k=0}^{\infty}z^k\vc{Y}(k)$, whose convergence radius is given
by $\sup\{y > 0; \sum_{k=0}^{\infty}y^k\vc{Y}(k) < \infty\}$. The
sequence $\{\vc{Y}(k)\}$ is said to be light-tailed if the convergence
radius of $\vc{Y}^{\ast}(z)$ is greater than one.
\end{defn}

\begin{rem}\label{rem-light-tailed}
Let $\overline{\vc{Y}}(k) = \sum_{l=k+1}^{\infty} \vc{Y}(l)$ for
$k=0,1,\dots$, and let $\overline{\vc{Y}}^{\ast}(z)$ denote the
generating function of $\{\overline{\vc{Y}}(k)\}$ defined by
$\sum_{k=0}^{\infty}z^k \overline{\vc{Y}}(k)$. We then have
\[
\overline{\vc{Y}}^{\ast}(z) 
= {\vc{Y}^{\ast}(1) - \vc{Y}^{\ast}(z) \over 1 - z},
\]
which implies that $\{\overline{\vc{Y}}(k)\}$ is light-tailed if and
only if $\{\vc{Y}(k)\}$ is light-tailed.
\end{rem}

Let $\overline{\vc{x}}^{\ast}(z)$ and $\vc{x}^{\ast}(z)$ denote the
generating functions of $\{\overline{\vc{x}}(k); k=0,1,\dots\}$ and
$\{\vc{x}(k); k=1,2,\dots\}$, respectively. We then have
\begin{equation}
\overline{\vc{x}}^{\ast}(z)
= {\vc{x}^{\ast}(1) - \vc{x}^{\ast}(z) \over 1 - z}.
\label{eqn-overline{x}^{ast}(z)}
\end{equation}
It also
follows from (\ref{eqn-x(k)}) that
\begin{equation}
\vc{x}^{\ast}(z)(\vc{I} - \vc{A}^{\ast}(z)/z) 
= (\vc{x}(0)\vc{B}^{\ast}(z) - \vc{x}(1)\vc{A}(0)),
\label{eqn-{x}^{ast}(z)}
\end{equation}
where $\vc{A}^{\ast}(z)$ and $\vc{B}^{\ast}(z)$ denote the generating
functions of $\{\vc{A}(k)\}$ and $\{\vc{B}(k)\}$ defined by
$\sum_{k=0}^{\infty} z^k \vc{A}(k)$ and $\sum_{k=1}^{\infty}z^k
\vc{B}(k)$, respectively.

Let $r_A$ and $r_B$ denote the convergence radii of $\vc{A}^{\ast}(z)$
and $\vc{B}^{\ast}(z)$, respectively. We then make the second
assumption.
\begin{assumption}\label{assu-2}
(a) $r_A > 1$, and (b) $r_B > 1$.
\end{assumption}

Proposition \ref{prop-light-tailed} below follows from Remark
\ref{rem-light-tailed}, Theorem 3.1 in~\cite{Li05a} and Theorem 2 in
\cite{Li05b}.
\begin{prop}\label{prop-light-tailed}
Under Assumption \ref{assu-1}, $\{\overline{\vc{x}}(k)\}$ is
light-tailed if and only if Assumption \ref{assu-2} holds.
\end{prop}

For further discussion, we introduce the following definition.
\begin{defn}\label{def-delta}
For any finite square matrix $\vc{X}$ with (possibly) complex number
elements, let $\delta(\vc{X})$ denote a maximum-modulus eigenvalue of
$\vc{X}$, whose argument $\arg \delta(\vc{X})$ is nonnegative, and
whose real part ${\rm Re}\, \delta(\vc{X})$ is not less than those of
the other eigenvalues of maximum modulus. Clearly, $|\delta(\vc{X})|$
is the spectral radius of $\vc{X}$. In addition, if $\vc{X}$ is
nonnegative and irreducible, $\delta(\vc{X}) > 0$ is the
Perron-Frobenius eigenvalue of $\vc{X}$.
\end{defn}

We now make the third assumption, Assumption~\ref{assu-3}
below. Unless otherwise stated, Assumptions~\ref{assu-1}, \ref{assu-2}
and \ref{assu-3} are valid throughout this paper. However we will not
assume unnecessary conditions in each of the propositions, lemmas,
theorems and corollaries presented in the rest of this paper.
\begin{assumption}\label{assu-3}
There exists some finite $\theta$ such that $1 < \theta < r_A$ and
$\theta = \delta(\vc{A}^{\ast}(\theta))$.
\end{assumption}

\begin{rem}\label{rem-theta}
It is easy to see that $\delta(\vc{A}^{\ast}(1))= 1$ and $\delta(\vc
{A}^{\ast}(e^s))$ $(s < \log r_A)$ is the Perron-Frobenius eigenvalue
of $\vc{A}^{\ast}(e^s)$. It is known \cite{King61} that
$\delta(\vc{A}^{\ast}(e^s))$ is a non-decreasing convex function of
$s$ (see also Proposition~7 in \cite{Gail98}).  Further it follows
from Lemma~2.3.3 in \cite{Neut89} that under
Assumption~\ref{assu-1}~(b), $\lim_{s\uparrow 0}(\rmd/\rmd
s)\delta(\vc{A}^{\ast}(e^s)) = \rho$. Therefore
Assumption~\ref{assu-1}~(b), (c) and Assumption~\ref{assu-2}~(a) are a
necessary condition for Assumption~\ref{assu-3}, though they are not
sufficient. In fact, $\lim_{y\uparrow
  r_A}\delta(\vc{A}^{\ast}(y))/y\le1$ in some cases (such an example
is given in Appendix \ref{appendix-theta}). A sufficient condition for
Assumption~\ref{assu-3} can be found in Theorem~4.12 in \cite{Bini05},
and see also Section~3 in~\cite{Gail98}. Finally, it should be noted
that Assumption~\ref{assu-3} requires
\begin{equation}
\left.{\rmd \over \rmd y} \delta(\vc{A}^{\ast}(y))
\right|_{y=\theta} - 1 > 0,
\label{add-ineqn-01}
\end{equation}
which will be used to show that the prefactors of the asymptotic
formulae presented in Section~\ref{sec-main-results} are positive.
\end{rem}

As is well known, the constant $\theta$ in Assumption \ref{assu-3}
plays a role in the light-tailed asymptotic analysis of
$\{\overline{\vc{x}}(k)\}$.  Several researchers have studied the
light-tailed asymptotics of $\{\overline{\vc{x}}(k)\}$ under the
assumption of $\theta < r_B$, where $\{\overline{\vc{x}}(k)\}$ decays
geometrically with rate $1/\theta$. Using the Tauberian theorem, Abate
et al. \cite{Abat94} presented a necessary condition for
\begin{equation}
\lim_{k\to\infty} \theta^k \vc{x}(k) = \vc{d},
\label{asymp-formula-01}
\end{equation}
where $\vc{d}$ is some positive vector. Note here that
(\ref{asymp-formula-01}) yields
\begin{equation}
\lim_{k\to\infty} \theta^k \overline{\vc{x}}(k) 
= (\theta - 1)^{-1}\vc{d}.
\label{asymp-formula-01'}
\end{equation}
M{\o}ller~\cite{Moll01} studied the asymptotic formula
(\ref{asymp-formula-01}) by considering the inter-visit times of level
zero, though his approach does not yield an explicit expression of
$\vc{d}$. Falkenberg~\cite{Falk94} and Gail et al.~\cite{Gail00}
obtained the asymptotic formula (\ref{asymp-formula-01}) by locating
the {\it dominant poles} (i.e., the maximum-order minimum-modulus
poles; see Definition~\ref{defn-dominant-pole}.)\ of the generating
function $\vc{x}^{\ast}(z)$ of $\{\vc{x}(k)\}$. However, Falkenberg's
sufficient condition for (\ref{asymp-formula-01}) includes a redundant
condition (see Remark~\ref{rem-Falk}).

Takine~\cite{Taki04} presented geometric asymptotic formulae of the
following form:
\begin{equation}
\lim_{n\to\infty} \theta^{nh+l} \overline{\vc{x}}(nh+l) = \vc{d}_l,
\qquad l=0,1,\dots,h-1,
\label{asymp-formula-02}
\end{equation}
where $h$ is some positive integer and $\vc{d}_l$'s
($l=0,1,\dots,h-1$) are some positive vectors. Clearly
(\ref{asymp-formula-02}) includes (\ref{asymp-formula-01'}) as a
special case. For simplicity, the cases of $h=1$ and $h\ge2$ are
called {\it the exactly geometric case} and {\it the periodically
  geometric case}, respectively.  Using the Markov renewal approach,
Takine~\cite{Taki04} derived the expression of $\vc{d}$ for the
exactly geometric case (i.e., $h=1$) and that of $\vc{d}_l$
($l=0,1,\dots,h-1$) for the periodically geometric case (i.e., $h \ge
2$), separately.

Li and Zhao~\cite{Li05b} studied the light-tailed asymptotics of the
stationary distribution of a Markov chain of GI/G/1 type. Corollary~1
therein would imply that the periodically geometric case is
impossible, which is inconsistent with Lemma~\ref{prop-smallest-pole}
in this paper. Therefore their results are valid only under the
condition that excludes the periodically geometric case.

In all the studies mentioned above, it is assumed that the phase space
is finite. Meanwhile, Miyazawa \cite{Miya04a}, Miyazawa and
Zhao~\cite{Miya04b} and Li et al.~\cite{LiHui07} considered structured
Markov chains with infinitely many phases, excluding the periodically
geometric case. Miyazawa \cite{Miya04a} and Miyazawa and
Zhao~\cite{Miya04b} studied Markov chains of M/G/1 type and GI/G/1
type, respectively, and they derived asymptotic formulae like
(\ref{asymp-formula-01}). Li et al.~\cite{LiHui07} derived an exactly
geometric asymptotic formula for the stationary distribution of a
quasi-birth-and-death (QBD) process and applied the obtained results
to a generalized join-the-shortest-queue model. The infiniteness of
the phase space causes some difficulties in asymptotic analysis of the
stationary tail probability vectors. That is a challenging problem,
but it is beyond the scope of this paper.

This paper studies the light-tailed asymptotics of
$\{\overline{\vc{x}}(k)\}$ of a Markov chain of M/G/1 type, not
excluding the periodically geometric case. The complex analysis
approach used in this paper is basically the same as Falkenberg
\cite{Falk94}'s approach, i.e., that is based on locating the dominant
poles of $\overline{\vc{x}}^{\ast}(z)$. Indeed, the approach is
classical, but it enables us to deal with the case of $\theta \ge r_B$
as well as the exactly and periodically geometric cases under the
condition $\theta < r_B$ in a unified manner. In addition, the complex
analysis approach gives us a deeper insight into the period in the
light-tailed asymptotics of $\{\overline{\vc{x}}(k)\}$, compared with
the Markov renewal approach \cite{Miya04a,Taki04}. In this paper, we
first present a simple and unified formula for both the exactly
geometric and periodically geometric cases, assuming $\theta <
r_B$. We also show that $h$ in (\ref{asymp-formula-02}) is closely
related to the period of a Markov additive process (MAdP) with kernel
$\{\vc{A}(k+1); k =0,\pm 1,\pm2,\dots\}$, where $\vc{A}(k) = \vc{O}$
for $k \le -1$. As for the case of $\theta \ge r_B$, we derive
light-tailed asymptotic formulae for $\{\overline{\vc{x}}(k)\}$ under
some mild conditions. We can find no previous studies paying special
attention to the case of $\theta \ge r_B$. As far as we know, only a
few examples of this case have been shown in Li and
Zhao~\cite{Li05b}. Therefore this paper is the first comprehensive
report of the light-tailed asymptotics in the case of $\theta \ge
r_B$.

The rest of this paper is organized as follows. In section
\ref{sec-preliminary}, we provide some preliminaries on
$\{\overline{\vc{x}}(k)\}$ and the period of a MAdP related to the
Markov chain of M/G/1 type. In section~\ref{sec-main-results}, we
present light-tailed asymptotic formulae for three cases: $\theta <
r_B$, $\theta > r_B$ and $\theta = r_B$ in subsections
\ref{subsec-case-1}, \ref{subsec-case-2} and \ref{subsec-case-3},
respectively.  Further in the appendix, we describe fundamental
results of the period of MAdPs, which play an important role in the
asymptotic analysis of the stationary distributions of structured
Markov chains such as ones of M/G/1 type and GI/G/1 type.

\section{Preliminaries}\label{sec-preliminary}

Throughout this paper, we use the following conventions. Let
$\Z=\{0,\pm1,\pm2,\dots\}$ and $\N = \{1,2,\dots\}$. Let $\C$ denote
the set of complex numbers. Let $\omega$ denote a complex number such
that $|\omega| =~1$. Let $\iota$ denote the imaginary unite, i.e.,
$\iota = \sqrt{-1}$. For any matrix $\vc{X}$ (resp.~vector $\vc{y}$),
its $(i,j)$th (resp.~$j$th) element is denoted by $[\vc{X}]_{i,j}$
(resp.~$[\vc{y}]_j$). When a (possibly complex-valued) function $f$
and a nonnegative function $g$ on $[0,\infty)$ satisfy $|f(x)| \le C
  g(x)$ for any sufficiently large $x$, we write $f(x) = O(g(x))$. We
  also write $f(x) = o(g(x))$ if $\lim_{x\to\infty}|f(x)|/g(x) = 0$.

\subsection{Some known results on the stationary tail probability vectors}
\label{subsec-M/G/1} 

Let $\vc{G}$ denote an $M \times M$ matrix whose $(i,j)$th ($i,j \in
\M$) element represents $\Pr[S_{a(k)} = j \mid X_0 = k+1, S_0 = i]$
for any given $k$ ($k\in\N)$, where $a(k) = \inf\{n \in \N; X_n =
k\}$. It is clear from the definition of $\vc{G}$ that
\begin{equation}
\vc{G}
=
\sum_{k=0}^{\infty}
\vc{A}(k)\vc{G}^k.
\label{eqn-G}
\end{equation}
It is known \cite{Neut89} that $\vc{G}$ is the minimal nonnegative
solution of $\vc{X} = \sum_{k=0}^{\infty} \vc{A}(k)\vc{X}^k$. If
Assumption~\ref{assu-1} (b) and (c) hold, $\vc{G}$ is stochastic (see
Theorem~2.3.1 in \cite{Neut89}).

The following result is an extension of Theorem~7.2.1 in \cite{Lato99}
to the Markov chain of M/G/1 type.
\begin{prop}\label{prop-structure-G}
If Assumption \ref{assu-1}~(a) and (b) hold, then $\vc{G}$ is
irreducible, or after some permutations it takes a form such that
\[
\vc{G} = 
\left(
\begin{array}{cc}
\vc{G}_1 & \vc{O}
\\
\vc{G}_{\bullet,1} & \vc{G}_{\bullet}
\end{array}
\right),
\]
where $\vc{G}_1$ is irreducible and $\vc{G}_{\bullet}$ is strictly
lower triangular.
\end{prop}

\proof See Appendix~\ref{proof-prop-structure-G}.\qed

\begin{rem}
The proof of Proposition~\ref{prop-structure-G} does not require that
$\vc{G}$ is stochastic. Thus Assumption~\ref{assu-1}~(c) is not
needed.
\end{rem}

Let $\vc{K}$ denote an $M_0 \times M_0$ matrix whose $(i,j)$th ($i,j
\in \M_0$) element represents 
$\Pr[S_{a(0)} = j \mid X_0 = 0, S_0 = i]$. Matrix $\vc{K}$ is given by
\[
\vc{K} = \vc{B}(0) + \sum_{k=1}^{\infty}\vc{B}(k)\vc{G}^{k-1}
\left(\vc{I} -
\sum_{m=1}^{\infty}\vc{A}(m)\vc{G}^{m-1}\right)^{-1}\vc{C}(0),
\]
and $\vc{x}(0)$ is given by 
\[
\vc{x}(0) = \left[
1 + {\vc{\kappa} \over 1 - \rho} 
\left\{ \vc{\beta}_B + 
\left(\vc{B} - \sum_{k=1}^{\infty}\vc{B}(k)\vc{G}^k \right)
[\vc{I} - \vc{A} + \vc{e}\vc{\pi}]^{-1} \vc{\beta}_A\right\}
\right]^{-1}\vc{\kappa},
\]
where $\vc{\kappa}$ denotes the stationary probability vector of
$\vc{K}$ (see Theorem~3.1 in~\cite{Sche90}). Further $\vc{x}(k)$
($k=1,2,\dots$) is determined by Ramaswami's \cite{Rama88}
recursion:
\begin{equation}
\vc{x}(k) = 
\left(
\vc{x}(0)\vc{U}_0(k) + \sum_{j=1}^{k-1}\vc{x}(j)\vc{U}(k-j)
\right)(\vc{I} - \vc{U}(0))^{-1},
\label{eqn-x(k)-G}
\end{equation}
where $\vc{U}(k)$ ($k=0,1,\dots$) and $\vc{U}_0(k)$ ($k=1,2,\dots$)
are given by
\begin{equation}
\vc{U}(k)
=
\sum_{m=k+1}^{\infty}
\vc{A}(m)\vc{G}^{m-k-1},
\qquad
\vc{U}_0(k)
=
\sum_{m=k}^{\infty}
\vc{B}(m)\vc{G}^{m-k},
\label{def-U-U_0}
\end{equation}
respectively. Note here that for any fixed $k \in \N$,
\begin{eqnarray*}
[\vc{U}(0)]_{i,j} 
&=& \Pr[S_{a(k)} = j, X_n \ge k~(n=1,2,\dots,a(k)) \mid X_0 = k, S_0=i],
\\ 
&&
\qquad \qquad \qquad \qquad \qquad \qquad \qquad 
\qquad \qquad \qquad ~~~i,j \in \M,
\end{eqnarray*}
and thus $\vc{I} - \vc{U}(0)$ is nonsingular due to Assumption
\ref{assu-1}~(a).

We now define $\vc{R}(k)$ and $\vc{R}_0(k)$ ($k \in \N$) as
\begin{equation}
\vc{R}(k)
=
\vc{U}(k)(\vc{I}-\vc{U}(0))^{-1},
\qquad
\vc{R}_0(k)
=
\vc{U}_0(k)(\vc{I}-\vc{U}(0))^{-1},
\label{def-R-R_0}
\end{equation}
respectively. For convenience, let $\vc{R}(0)=\vc{O}$. We then
rewrite (\ref{eqn-x(k)-G}) as
\begin{equation}
\vc{x}(k) = \vc{x}(0)\vc{R}_0(k)
+ \sum_{j=1}^k\vc{x}(j)\vc{R}(k-j).
\label{eqn-x(k)-R}
\end{equation}
Let $\vc{R}^{\ast}(z)$ and $\vc{R}_0^{\ast}(z)$ denote the generating
functions of $\{\vc{R}(k)\}$ and $\{\vc{R}_0(k)\}$ defined by
$\sum_{k=0}^{\infty}z^k\vc{R}(k)$ and
$\sum_{k=1}^{\infty}z^k\vc{R}_0(k)$, respectively. It then follows
from (\ref{eqn-x(k)-R}) that
\[
\vc{x}^{\ast}(z)
=
\vc{x}(0) \vc{R}_0^{\ast}(z)[\vc{I} - \vc{R}^{\ast}(z)]^{-1},
\]
from which and (\ref{eqn-overline{x}^{ast}(z)}) we have
\[
\overline{\vc{x}}^{\ast}(z)
= {\vc{x}^{\ast}(1) - \vc{x}(0) \vc{R}_0^{\ast}(z)
[\vc{I} - \vc{R}^{\ast}(z)]^{-1} \over 1 - z}.
\]
Using (\ref{eqn-G}), (\ref{def-U-U_0}) and
(\ref{def-R-R_0}), we can readily have the following result.
\begin{prop}\label{lemma-F^{ast}(z)}
If Assumption \ref{assu-1}~(a) and (b) hold, then
\begin{align}
\vc{I} - \vc{\Gamma}_A^{\ast}(z)
&=
(\vc{I} - \vc{R}^{\ast}(z))
(\vc{I}-\vc{U}(0))
(\vc{I}-\vc{G}/z),
&
0 &< |z| < r_A,
\label{RG-factorization}
\\
\vc{B}^{\ast}(z) - \vc{U}_0(1)\vc{G}
&= \vc{R}_0^{\ast}(z)(\vc{I}-\vc{U}(0))
(\vc{I}-\vc{G}/z),
& 0 &< |z| < r_B,
\label{eqn-R_0}
\end{align}
 where
\begin{equation}
\vc{\Gamma}_A^{\ast}(z) = z^{-1} \vc{A}^{\ast}(z).
\label{add-eqn-Gamma_A^{ast}(z)}
\end{equation}

\end{prop}

\begin{rem}
Equation (\ref{RG-factorization}) shows the $RG$-factorization of the
Markov chain of M/G/1 type (see, Theorem~14 in~\cite{Zhao03}).
\end{rem}

\begin{prop}\label{remark-RG-Fact}
Suppose Assumption~\ref{assu-1}~(a), (b) and
Assumption~\ref{assu-2}~(a) hold. Then
\begin{equation}
\det(\vc{I} - \vc{\Gamma}_A^{\ast}(z)) = 0~\mbox{if and only if}~
\det(\vc{I} - \vc{R}^{\ast}(z)) = 0,
\qquad 1 < |z| < r_A.
\label{eqn-RG}
\end{equation}
\end{prop}

\proof
Since $\vc{I} - \vc{U}(0)$ and $\vc{I} - \vc{G}/z$ are nonsingular for
$|z| > 1$, (\ref{RG-factorization}) leads to (\ref{eqn-RG}).
\qed

\medskip

It follows from (\ref{add-eqn-Gamma_A^{ast}(z)}) that
$\vc{\Gamma}_A^{\ast}(\theta) = \vc{A}^{\ast}(\theta)/\theta$ and thus
Assumption \ref{assu-3} yields
\begin{equation}
\delta(\vc{\Gamma}_A^{\ast}(\theta)) 
= \delta(\vc{A}^{\ast}(\theta))/\theta = 1.
\label{add-eqn-48}
\end{equation}
The matrix $\vc{\Gamma}_A^{\ast}(\theta)$ is irreducible, because
$[\vc{A}^{\ast}(\theta)]_{i,j} > 0$ if and only if
$[\vc{A}^{\ast}(1)]_{i,j} = [\vc{A}]_{i,j} > 0$. Therefore
$\delta(\vc{\Gamma}_A^{\ast}(\theta)) = 1$ is the Perron-Frobenius
eigenvalue of $\vc{\Gamma}_A^{\ast}(\theta)$. Let $\vc{\mu}(\theta)$
and $\vc{v}(\theta)$ denote the Perron-Frobenius left- and
right-eigenvectors of $\vc{\Gamma}_A^{\ast}(\theta)$, which are
normalized such that
\[
\vc{\mu}(\theta)\vc{e} = 1, 
\qquad 
\vc{\mu}(\theta)\vc{v}(\theta) = 1.
\]
Clearly $\vc{\mu}(\theta) > \vc{0}$ and $\vc{v}(\theta) >
\vc{0}$. Further (\ref{RG-factorization}) yields
$\vc{\mu}(\theta)\vc{R}^{\ast}(\theta) = \vc{\mu}(\theta)$. It thus
follows from Corollary~8.1.30 in \cite{Horn90} that
\begin{equation}
\delta(\vc{R}^{\ast}(\theta)) = 1.
\label{eqn-delta(R^*(theta))}
\end{equation}
Using this, we can prove the following result. 
\begin{prop}\label{prop-structure-R}
If Assumption \ref{assu-1}~(a)--(c), Assumption~\ref{assu-2}~(a) and
Assumption \ref{assu-3} hold, then $\vc{R} \triangleq
\vc{R}^{\ast}(1)$ is irreducible, or after some permutations it takes
a form such that
\[
\vc{R} = 
\left(
\begin{array}{cc}
\vc{R}_1 & \vc{R}_{1,\bullet}
\\
\vc{O} & \vc{R}_{\bullet}
\end{array}
\right),
\]
where $\vc{R}_1$ is irreducible and $\vc{R}_{\bullet}$ is strictly
upper triangular.
\end{prop}

\proof This proposition can be proved in a similar way to the proof of
Theorem~7.2.2 in \cite{Lato99}, but some modifications are needed. A
complete proof is given in Appendix \ref{proof-prop-G-R}.\qed

\begin{rem}
Because of some dual properties between $\vc{G}$ and $\vc{R}$, it
might seem that the statement of Proposition~\ref{prop-structure-R}
would be true under the conditions of
Proposition~\ref{prop-structure-G}, i.e., Assumption~\ref{assu-1}~(a)
and (b). In fact, if $\vc{T}$ is the transition probability matrix of
a QBD, the statement of Proposition~\ref{prop-structure-R} is true
under Assumption~\ref{assu-1}~(a) and (b) (see Theorems~7.2.1 and
7.2.2 in \cite{Lato99}). However, this is not the case in general.
Assume that $\vc{A}(k) = \vc{O}$ for all $k=2,3,\dots$ and the other
block matrices of $\vc{T}$ in (\ref{eqn-matrix-T}) are all positive.
Then Assumption~\ref{assu-1}~(a) and (b) are satisfied, whereas
$\vc{R} = \vc{O}$.
\end{rem}

Since $[\vc{R}]_{i,j} > 0$ implies $[\vc{R}^{\ast}(\theta)]_{i,j} > 0$
(vice versa), it follows from (\ref{eqn-delta(R^*(theta))}) and
Proposition \ref{prop-structure-R} that
$\delta(\vc{R}^{\ast}(\theta))=1$ is a simple eigenvalue.  Let
$\vc{s}(\theta) = (\vc{I}-\vc{U}(0))
(\vc{I}-\vc{G}/\theta)\vc{v}(\theta)$.  From (\ref{RG-factorization}),
we then have $\vc{R}^{\ast}(\theta)\vc{s}(\theta) = \vc{s}(\theta)$.
Note here that
\begin{eqnarray*}
(\vc{I}-\vc{G}/\theta)^{-1}(\vc{I}-\vc{U}(0))^{-1}\vc{s}(\theta)
&=& \vc{v}(\theta) > \vc{0},
\\
(\vc{I}-\vc{G}/\theta)^{-1}(\vc{I}-\vc{U}(0))^{-1} \ge \vc{O}
&\neq& \vc{O}.
\end{eqnarray*}
Thus $\vc{s}(\theta) \ge \vc{0},\neq \vc{0}$ (see, e.g., Theorem~8.3.1
in \cite{Horn90}). The following proposition summarizes the results on
the spectral radius of $\vc{R}^{\ast}(\theta)$ and its corresponding
eigenvectors.
\begin{prop}[Lemma 5 in \cite{Li05b}]\label{lemma5-Li05b}

Suppose Assumption \ref{assu-1}~(a)--(c), Assumption~\ref{assu-2}~(a)
and Assumption \ref{assu-3} hold.  Then $\vc{\mu}(\theta) > \vc{0}$
and $\vc{s}(\theta) \ge \vc{0},\neq \vc{0}$, which are left- and
right-eigenvectors, respectively, of $\vc{R}^{\ast}(\theta)$
corresponding simple eigenvalue $\delta(\vc{R}^{\ast}(\theta))=1$.
\end{prop}

We define $\vc{\Gamma}_R(k)$ ($k\in\Z$) as
\begin{equation}
\vc{\Gamma}_R(k)
= \left\{
\begin{array}{ll}
\theta^k
\diag(\vc{\mu}(\theta))^{-1}\vc{R}(k)^{\rm t}\diag(\vc{\mu}(\theta)), 
& k=1,2,\dots,
\\
\vc{O}, & k=0,-1,-2,\dots,
\end{array}
\right.
\label{def-Gamma_R(k)}
\end{equation}
where super-subscript ${\rm t}$ denotes transpose and
$\diag(\vc{\mu}(\theta))$ denotes an $M \times M$ diagonal matrix
whose $j$th diagonal element is equal to $[\vc{\mu}(\theta)]_j$. Let
$\vc{\Gamma}_R^{\ast}(z)$ denote the generating function of
$\{\vc{\Gamma}_R(k)\}$ defined by $\sum_{k \in \Z} z^k
\vc{\Gamma}_R(k)$. We then have
\begin{equation}
\vc{\Gamma}_R^{\ast}(z)
=
\diag(\vc{\mu}(\theta))^{-1} \vc{R}^{\ast}(\theta z)^{\rm t}
\diag(\vc{\mu}(\theta)).
\label{eqn-Gamma_R^*(z)}
\end{equation}
It is easy to see that $\vc{\Gamma}_R^{\ast}(1)$ is stochastic.
Further it follows from Propositions~\ref{prop-structure-R} and
\ref{lemma5-Li05b} that either $\vc{\Gamma}_R^{\ast}(1)$ is
irreducible or after some permutations, $\vc{\Gamma}_R^{\ast}(z)$
takes a form such that
\[
\vc{\Gamma}_R^{\ast}(z)
= \bordermatrix{
& \M^{(1)} 
& \M^{(2)} 
\cr
\M^{(1)} 
& \vc{\Gamma}_{R,1}^{\ast}(z) 
& \vc{O}
\cr
\M^{(2)}
& \vc{\Gamma}_{R,2,1}^{\ast}(z) 
& \vc{\Gamma}_{R,2}^{\ast}(z)
},
\]
where $\vc{\Gamma}_{R,1}^{\ast}(1)$ is irreducible and stochastic and
$\vc{\Gamma}_{R,2}^{\ast}(z)$ is strictly lower triangular, i.e.,
$\delta(\vc{\Gamma}_{R,2}^{\ast}(z)) = 0$. Thus according to
Theorem~\ref{append-theorem}, we define $h$ as
\begin{equation}
h 
= \max\{n\in \N; 
\delta(\vc{\Gamma}_R^{\ast}(\omega_n)) = 1\},
\label{eqn-h-01}
\end{equation}
where $\omega_x = \exp(2\pi\iota/x)$ for $x \ge 1$.  Note that $h$ is
the period of the recurrent class $\M^{(1)}$ in an MAdP
$\{\vc{\Gamma}_R(k);k \in \Z\}$. Let $h_{i,j}$ ($i,j \in \M$) denote
the first jump point of the distribution of the first passage time
from state $i$ to state $j$ in the MAdP $\{\vc{\Gamma}_R(k);k \in
\Z\}$. Then Takine~\cite{Taki04}'s asymptotic formulae are as follows.

\begin{prop}[Theorem 2 in~\cite{Taki04}]\label{prop-takine}

Suppose Assumptions \ref{assu-1}--\ref{assu-3} hold and $\theta <
r_B$. Let $\vc{\pi}_{\ast} = \sum_{k=1}^{\infty}\vc{x}(k)$, which is
given by (see Lemma~3 in~\cite{Taki04})
\[
\vc{\pi}_{\ast} 
= [\vc{x}(0)\{ \vc{B} + \vc{\beta}_B\vc{g} \} - \vc{x}(1)\vc{A}(0)] 
(\vc{I} - \vc{A} + (\vc{e} - \vc{\beta}_A)\vc{g})^{-1},
\]
where $\vc{g}$ denotes the stationary probability vector of
$\vc{G}$. If $h\ge2$, then for $l=0,1,\dots,h-1$,
\[
\lim_{n\to\infty}\theta^{nh+l}[\overline{\vc{x}}(nh+l)]_j 
= [\vc{d}_l]_j,
\qquad j \in \M,
\]
where 
\begin{eqnarray*}
[\vc{d}_l]_j
&=& {h \over (\rmd/\rmd y)\delta(\vc{A}^{\ast}(y))|_{y=\theta} - 1}
\sum_{\nu \in \M} \sum_{k=0}^{\infty} \theta^{l+kh-h_{j,\nu}} 
\nonumber
\\
&& {} \hspace{-3mm}
\times \sum_{m=l+kh-h_{j,\nu}+1}^{\infty} \hspace{-3mm}
 [\vc{x}(0)\vc{R}_0(m) + \vc{\pi}_{\ast} \vc{R}(m)]_{\nu} 
 [(\vc{I} - \vc{U}(0)) 
(\vc{I} - \vc{G}/\theta)\vc{v}(\theta)]_{\nu} [\vc{\mu}(\theta)]_j.
\end{eqnarray*}
\end{prop}

\begin{prop}[Theorem 3 in~\cite{Taki04}]
\label{prop-takine-2}

Suppose Assumptions \ref{assu-1}--\ref{assu-3} hold and $\theta <
r_B$. If $h = 1$, then
\begin{equation}
\lim_{k\to\infty}\theta^{k}\overline{\vc{x}}(k) 
= 
\dm{ \left[
\vc{x} (0) \vc{B}^{\ast}(\theta) - \vc{x}(1) \vc{A}(0) 
\right] \vc{v}(\theta)
\over (\theta - 1)
\{(\rmd/\rmd y)\delta(\vc{A}^{\ast}(y))|_{y=\theta} - 1\}
}
\cdot \vc{\mu}(\theta).
\label{add-eqn-50}
\end{equation}
\end{prop}

\subsection{Period of a related Markov additive process}
\label{subsec-MAdP}

 We consider a MAdP $\{(\breve{X}_n, \breve{S}_n);n=0,1,\dots\}$ with
 state space $\Z \times \M$ and kernel $\{\vc{\Gamma}_A(k);k\in \Z\}$,
 where
\begin{equation}
\vc{\Gamma}_A(k)
=
\left\{
\begin{array}{ll}
\vc{A}(k+1), & k=-1,0,1,\dots,
\\
\vc{O}, & k=-2,-3,-4,\dots.
\end{array}
\right.
\label{def-Gamma_A(k)}
\end{equation}
It follows from (\ref{add-eqn-Gamma_A^{ast}(z)}) and
(\ref{def-Gamma_A(k)}) that $\sum_{k \in \Z}z^k \vc{\Gamma}_A(k) =
\vc{\Gamma}_A^{\ast}(z)$. It is easy to see that for $i,j \in \M$,
\begin{eqnarray}
&& \Pr[\breve{X}_{n+1}=k_{\ast}+k, \breve{S}_{n+1} = j 
\mid \breve{X}_n=k_{\ast}, \breve{S}_n = i]
\nonumber
\\
&& \quad {} 
= \Pr[X_{n+1}=k_{\ast}+k, S_{n+1} = j \mid X_n=k_{\ast}, S_n = i],
\quad k_{\ast}+k \in \N,~k_{\ast} \in \N.
\nonumber
\\
\label{add-eqn-58}
\end{eqnarray}
For any two states $(k_1,j_1)$ and $(k_2,j_2)$ in $\Z \times \M$, we
write $(k_1,j_1) \rightarrow (k_2,j_2)$ when there exists a path from
$(k_1,j_1)$ to $(k_2,j_2)$ with some positive probability.
\begin{prop}\label{prop-period}
Suppose Assumption \ref{assu-1} (a) and (b) hold. Then for each $j \in
\M$, there exists a nonzero integer $k_j$ such that $(0,j) \rightarrow
(k_j,j)$.
\end{prop}

\proof See Appendix
\ref{proof-lem-period}. \qed

\medskip

Assumption \ref{assu-1} (b) and (\ref{add-eqn-Gamma_A^{ast}(z)}) show
that $\vc{\Gamma}_A^{\ast}(1) = \vc{A}$ is irreducible, from which and
Proposition~\ref{prop-period} it follows that the period of the MAdP
$\{(\breve{X}_n, \breve{S}_n);n=0,1,\dots\}$ is well-defined and is
denoted by $\tau$ (see Definition~\ref{def-period-MAdP-01}).  Note
here that from (\ref{add-eqn-58}) and the definition of $\vc{G}$,
\[
[\vc{G}]_{i,j} = \Pr[\breve{S}_{\breve{a}(k)} = j \mid
  \breve{X}_0=k+1, \breve{S}_0=i],
\qquad
i,j \in \M,
\]
where $\breve{a}(k) = \inf\{n \in \N; \breve{X}_n = k\}$. Note also
that $\vc{G}$ has no zero rows due to
Assumption~\ref{assu-1}~(a). Proposition~\ref{prop-structure-G} then
implies that there exists a
path such that
\[
(k,j_0) 
\rightarrow (k-1,j_1)
\rightarrow (k-2,j_2)
\rightarrow \cdots \rightarrow 
(k-M,j_M),
\]
where $k \in \Z$ and $j_n \in \M$ for $n=0,1,\dots,M$. In the above
path, a phase appears at least two times, and thus $\tau \le M$.

\begin{prop}[Propositions~13 and 14 in \cite{Gail96}]\label{prop-tau}
Under Assumption \ref{assu-1} (a) and (b),
\[
\tau = \max\{n\in\M; \delta(\vc{\Gamma}_A^{\ast}(\omega_n)) = 1\}.
\]
\end{prop}

\medskip

Lemma~\ref{def-period-MAdP-02} shows that there exists a function\footnote{corrected: ``an injective function" $\longrightarrow$ ``a function"} $p$ from $\M$ to $\{0,1,\dots,\tau-1\}$ such that
\[
[\vc{\Gamma}_A(k)]_{i,j} > 0
\mbox{ only if } k \equiv
p(j)-p(i)~(\Mod~\tau).
\]
Let $\vc{\Delta}_M(z)$ denote an $M \times M$ diagonal matrix such
that
\[
\vc{\Delta}_M(z)
= 
\left(
\begin{array}{cccc}
z^{-p(1)} &          &         &
\\
         & z^{-p(2)} &         &
\\
         &          &\ddots   &  
\\
         &          &         & z^{-p(M)}
\end{array}
\right).
\]
Let $\vc{\mu}(z)$ and $\vc{v}(z)$ ($z \in \C$) denote left- and
right-eigenvectors of $\vc{\Gamma}_A^{\ast}(z)$ corresponding to
eigenvalue $\delta(\vc{\Gamma}_A^{\ast}(z))$, normalized such that
\[
\vc{\mu}(z) \vc{\Delta}_M(z/|z|)\vc{e} = 1, 
\qquad 
\vc{\mu}(z)\vc{v}(z) = 1.
\]
Note that if $z$ is real, $\vc{\Delta}_M(z/|z|) = \vc{I}$ and
therefore the definition of $\vc{\mu}(z)$ and $\vc{v}(z)$ is
consistent with that of $\vc{\mu}(\theta)$ and $\vc{v}(\theta)$.

The following proposition is an immediate
consequence of Lemma \ref{lem-delta-01}.
\begin{prop}\label{prop-delta-Gamma_A-1}
If Assumption \ref{assu-1} (a) and (b) hold, then for $0 < y < r_A$,
\begin{align}
&&
\delta(\vc{\Gamma}_A^{\ast}(y\omega_{\tau}^{\nu})) 
&= \delta(\vc{\Gamma}_A^{\ast}(y)),
& \nu=0,1,\dots,\tau-1,
&&
\label{add-eqn-35}
\\
&&
\vc{\mu}(y\omega_{\tau}^{\nu})
&= \vc{\mu}(y)\vc{\Delta}_M(\omega_{\tau}^{\nu})^{-1},
& \nu=0,1,\dots,\tau-1,
&&
\label{add-eqn-36a}
\\
&&
\vc{v}(y\omega_{\tau}^{\nu}) 
&= \vc{\Delta}_M(\omega_{\tau}^{\nu})\vc{v}(y),
& \nu=0,1,\dots,\tau-1.
&&
\label{add-eqn-36b}
\end{align}

\end{prop}

It follows from (\ref{add-eqn-Gamma_A^{ast}(z)}) that
\begin{equation}
\delta(\vc{\Gamma}_A^{\ast}(y)) = y^{-1}\delta(\vc{A}^{\ast}(y)),
\qquad 0 <  y < r_A.
\label{add-eqn-53}
\end{equation}
Thus according to the property of $\delta(\vc{A}^{\ast}(y))$ (see
Remark~\ref{rem-theta}), we obtain
\begin{eqnarray}
\delta(\vc{\Gamma}_A^{\ast}(y)) &<& 1, 
\qquad 1 < y < \theta.
\label{add-eqn-20}
\end{eqnarray}
Further from (\ref{add-eqn-53}), we have
\begin{equation}
\left.{\rmd \over \rmd z} \delta(\vc{\Gamma}_A^{\ast}(z))
\right|_{z=\theta}
= \left.{\rmd \over \rmd y} \delta(\vc{\Gamma}_A^{\ast}(y))
\right|_{y=\theta}
= {1 \over \theta} \left(\left.{\rmd \over \rmd y} \delta(\vc{A}^{\ast}(y))
\right|_{y=\theta} - 1\right) > 0,
\label{eqn-deff-delta}
\end{equation}
where the second equality follows from (\ref{add-eqn-48}) and the last
inequality follows from (\ref{add-ineqn-01}).

The following proposition can be easily obtained by
(\ref{add-eqn-48}), Proposition~\ref{prop-tau} and Theorem
\ref{append-theorem}.
\begin{prop}\label{prop-delta-Gamma_A-2}
Suppose Assumption \ref{assu-1}~(a)--(c), Assumption \ref{assu-2}~(a)
and Assumption~\ref{assu-3} hold. Then
$\delta(\vc{\Gamma}_A^{\ast}(\theta\omega)) = 1$ if and only if
$\omega^{\tau} = 1$. Thus
\[
\tau = \max\{n \in \M;
\delta(\vc{\Gamma}_A^{\ast}(\theta\omega_n))=1\}.
\]
Further if $\delta(\vc{\Gamma}_A^{\ast}(\theta\omega)) = 1$, the
eigenvalue is simple.
\end{prop}

\begin{prop}\label{prop-h-tau}
If Assumption \ref{assu-1}~(a)--(c), Assumption \ref{assu-2}~(a) and
Assumption~\ref{assu-3} hold, then $\tau$ is equal to period $h$ of
the MAdP with kernel $\{\vc{\Gamma}_{R,1}(k); k\in\Z\}$.
\end{prop}

\proof It follows from (\ref{eqn-h-01}) and
Proposition~\ref{prop-delta-Gamma_A-2} that $h = \tau$ if the
following is true.
\begin{equation}
\delta(\vc{\Gamma}_R^{\ast}(\omega)) = 1
~\mbox{if and only if}~
\delta(\vc{\Gamma}_A^{\ast}(\theta\omega)) = 1.
\label{add-eqn-45}
\end{equation}
In fact, (\ref{eqn-Gamma_R^*(z)}) shows that
\begin{equation}
\delta(\vc{\Gamma}_R^{\ast}(\omega)) = 1~\mbox{if and only if}~
\delta(\vc{R}^{\ast}(\theta\omega)) = 1.
\label{add-eqn-43}
\end{equation}
Since $\delta(\vc{\Gamma}_A^{\ast}(\theta))
=\delta(\vc{R}^{\ast}(\theta)) = 1$,
$|\delta(\vc{\Gamma}_A^{\ast}(\theta\omega))| \le 1$ and
$|\delta(\vc{R}^{\ast}(\theta\omega))| \le 1$ (see Theorem~8.1.18 in
\cite{Horn90}). Therefore Proposition \ref{remark-RG-Fact} implies
that
\[
\delta(\vc{R}^{\ast}(\theta\omega)) = 1
~\mbox{if and only if}~\delta(\vc{\Gamma}_A^{\ast}(\theta\omega)) = 1.
\]
This and (\ref{add-eqn-43}) yield (\ref{add-eqn-45}).
\qed


\section{Main Results}\label{sec-main-results}

Proposition \ref{prop-light-tailed} and Definition
\ref{def-light-tailed} show that $[\overline{\vc{x}}^{\ast}(z)]_j$ ($j
\in \M$) is holomorphic in the domain $\{z\in\C; |z| \le 1\}$ and thus
has no pole in the same domain. It follows from
(\ref{eqn-overline{x}^{ast}(z)}), (\ref{eqn-{x}^{ast}(z)}) and
(\ref{add-eqn-Gamma_A^{ast}(z)}) that for $z$ such that $[ \vc{I} -
  \vc{\Gamma}_A^{\ast}(z) ]^{-1}$ exists,
\begin{equation}
\overline{\vc{x}}^{\ast}(z) 
=
{\vc{x}^{\ast}(1) \over 1 - z} 
-
{\vc{x} (0) \vc{B}^{\ast}(z) - \vc{x}(1) \vc{A}(0)  \over 1 - z}
\left[ \vc{I} - \vc{\Gamma}_A^{\ast}(z) \right]^{-1}.
\label{add-eqn-x^{ast}(z)-02}
\end{equation}
By definition,
\begin{equation}
\left[ \vc{I} - \vc{\Gamma}_A^{\ast}(z) \right]^{-1}
= {\adj(\vc{I} - \vc{\Gamma}_A^{\ast}(z)) 
\over \det(\vc{I} - \vc{\Gamma}_A^{\ast}(z))},
\label{eqn-inverse}
\end{equation}
where $\adj(\vc{I} - \vc{\Gamma}_A^{\ast}(z))$ denotes the adjugate of
$\vc{I} - \vc{\Gamma}_A^{\ast}(z)$.  Therefore the roots of
$\det(\vc{I} - \vc{\Gamma}_A^{\ast}(z)) = 0$ with $|z|>1$, if any, are
candidates for the dominant poles of $[\overline{\vc{x}}^{\ast}(z)]_j$
($j \in \M$).

Let $r_m(z)$'s ($m=1,2,\dots,M$)
denote the eigenvalues of $\vc{\Gamma}_A^{\ast}(z)$ such that
\begin{eqnarray}
&& r_1(z) = \delta(\vc{\Gamma}_A^{\ast}(z)),
\qquad
|r_1(z)| \ge |r_2(z)| \ge \cdots \ge |r_M(z)|,
\\
&& \det(\vc{I} - \vc{\Gamma}_A^{\ast}(z)) 
= \prod_{m=1}^M (1 - r_m(z)).
\label{eqn-det}
\label{def-r_1(z)}
\end{eqnarray}
For convenience in what follows, let $\varepsilon_0$
denote a sufficiently small positive number, which may take different
values in different places.
\begin{lem}\label{prop-smallest-pole}
Suppose Assumption \ref{assu-1}~(a)--(c), Assumption \ref{assu-2}~(a)
and Assumption~\ref{assu-3} hold. Then the equation $\det(\vc{I} -
\vc{\Gamma}_A^{\ast}(z)) = 0$ has exactly $\tau$ roots
$\theta\omega_{\tau}^{\nu}$ ($\nu=0,1,\dots,\tau-1$) in the domain
$\{z\in\C;1 < |z| < \theta+\varepsilon_0\}$. In addition, each of the
roots is simple.
\end{lem}

\proof It follows from (\ref{add-eqn-35}) and (\ref{eqn-deff-delta})
that $(\rmd/\rmd z)\delta(\vc{\Gamma}_A^{\ast}(z))
|_{z=\theta\omega_{\tau}^{\nu}} > 0$ for $\nu=0,1,\dots,\tau-1$. Thus,
according to the former part of
Proposition~\ref{prop-delta-Gamma_A-2},
$\{\theta\omega_{\tau}^{\nu};\nu=0,1,\dots,\tau-1\}$ are simple roots
of the equation $\delta(\vc{\Gamma}_A^{\ast}(z)) - 1 = 0$. Further the
latter part of the proposition implies that $\prod_{m=2}^M (1 -
r_m(\theta\omega_{\tau}^{\nu})) \neq 0$ for
$\nu=0,1,\dots,\tau-1$. Therefore
$\{\theta\omega_{\tau}^{\nu};\nu=0,1,\dots,\tau-1\}$ are simple roots
of the equation $\det(\vc{I} - \vc{\Gamma}_A^{\ast}(z)) = 0$. Note
here that $\det(\vc{I} - \vc{\Gamma}_A^{\ast}(z))$ is holomorphic and
not identically zero in the domain $\{z\in\C;0 < |z| < r_A\}$. As a
result, it suffices to show that $\det(\vc{I} -
\vc{\Gamma}_A^{\ast}(z))=0$ has no other roots in the domain
$\{z\in\C;1 < |z| \le \theta\}$.

By definition, if $\delta(\vc{\Gamma}_A^{\ast}(\theta\omega)) \neq 1$,
then $r_m(\theta\omega) \neq 1$ for all $m=2,3,\dots,M$. Thus
$\det(\vc{I} - \vc{\Gamma}_A^{\ast}(\theta\omega)) = 0$ only if
$\omega = \omega_{\tau}^{\nu}$ ($\nu=0,1,\dots,\tau-1$).  In addition,
it follows from (\ref{add-eqn-20}) and Theorem~8.1.18 in~\cite{Horn90}
that for $m=2,3,\dots,M$,
\[
|r_m(z)|
\le |\delta(\vc{\Gamma}_A^{\ast}(z))| 
\le \delta(\vc{\Gamma}_A^{\ast}(|z|)) < 1, \qquad 
1 < |z| < \theta.
\]
This and (\ref{def-r_1(z)}) show that $\det(\vc{I} -
\vc{\Gamma}_A^{\ast}(z)) \neq 0$ in the domain $\{z\in\C;1 < |z| <
\theta\}$.  \qed

\medskip

\begin{rem}
In the proof of Corollary~1 in \cite{Li05b}, it is stated that
$\det(\vc{I} - \vc{\Gamma}_A^{\ast}(z)) = 0$ has one and only one root
on the circle $\{z\in\C;|z|=\theta\}$, which is, in general,
incorrect.
\end{rem}

\begin{lem}\label{lem-adj}
Suppose Assumptions \ref{assu-1}~(a)--(c), Assumption~\ref{assu-2}~(a)
and Assumption~\ref{assu-3} hold. Then for $\nu=0,1,\dots,\tau-1$,
\begin{eqnarray}
\adj( \vc{I} - \vc{\Gamma}_A^{\ast}(\theta \omega_{\tau}^{\nu}))
=
\prod_{m=2}^{M}(1 - r_m(\theta \omega_{\tau}^{\nu})) 
\vc{v}(\theta \omega_{\tau}^{\nu}) \vc{\mu}(\theta \omega_{\tau}^{\nu})
\neq \vc{O},
\nonumber
\\
\qquad \nu = 0,1,\dots,\tau-1.
\label{eq-adj-06}
\end{eqnarray}
\end{lem}

\proof See Appendix \ref{proof-lem-adj}. \qed

\medskip

In the rest of this section, we first derive a light-tailed asymptotic
formula for the case of $\theta < r_B$ in
subsection~\ref{subsec-case-1}. We then discuss the cases of $\theta >
r_B$ and $\theta = r_B$ in subsections~\ref{subsec-case-2} and
\ref{subsec-case-3}, respectively.

\subsection{Case of $\theta < r_B$}\label{subsec-case-1}

Lemmas~\ref{prop-smallest-pole} and \ref{lem-adj} imply that if
$\theta < \min(r_A,r_B)$, $[\overline{\vc{x}}^{\ast}(z)]_j$ ($j\in\M$)
is holomorphic for $|z| < \theta$ and meromorphic for $|z| \le
\theta$, and thus the candidates for the dominant poles of
$[\overline{\vc{x}}^{\ast}(z)]_j$ are the simple roots
$\{\theta\omega_{\tau}^{\nu};\nu=0,1,\dots,\tau-1\}$ of $\det(\vc{I}
- \vc{\Gamma}_A(z)) = 0$. Note here that $z=\theta\omega_{\tau}^{\nu}$
($\nu=0,1,\dots,\tau-1$) is a simple pole of
$[\overline{\vc{x}}^{\ast}(z)]_j$ ($j\in\M$) if and only if
\[
0 < \lim_{z\to\theta\omega_{\tau}^{\nu}} 
\left|
\left(1 - {z\over \theta\omega_{\tau}^{\nu}} \right)
[\overline{\vc{x}}^{\ast}(z)]_j
\right| < \infty.
\]
Thus it follows from Theorem~\ref{append-thm-asymp} and
Remark~\ref{rem-theorem-A.1} that if each
$[\overline{\vc{x}}^{\ast}(z)]_j$ ($j \in \M$) has at least one pole
of $\{\theta\omega_{\tau}^{\nu};\nu=0,1,\dots,\tau-1\}$, then
\begin{equation}
\overline{\vc{x}}(k)
= \theta^{-k}\sum_{\nu=0}^{\tau-1} {1 \over (\omega_{\tau}^{\nu})^k}
\lim_{z\to\theta\omega_{\tau}^{\nu}} 
\left(1 - {z\over \theta\omega_{\tau}^{\nu}} \right)
\overline{\vc{x}}^{\ast}(z)
+ O((\theta+\varepsilon_0)^{-k})\vc{e}^{\rm t},
\label{add-eqn-40}
\end{equation}
where the dominant term on the right hand side of (\ref{add-eqn-40})
is positive for any $k=0,1,\dots$ (see
Theorem~\ref{append-thm-asymp}~(d)), i.e., for any $j \in \M$ and
$k=0,1,\dots$,
\begin{equation}
\left[
\theta^{-k}\sum_{\nu=0}^{\tau-1} {1 \over (\omega_{\tau}^{\nu})^k}
\lim_{z\to\theta\omega_{\tau}^{\nu}} 
\left(1 - {z\over \theta\omega_{\tau}^{\nu}} \right)
\overline{\vc{x}}^{\ast}(z)
\right]_j > 0.
\label{add-eqn-57}
\end{equation}

We now note that (\ref{add-eqn-x^{ast}(z)-02}) yields
\begin{eqnarray}
\lim_{z\to\theta\omega_{\tau}^{\nu}} 
\left(1 - {z\over \theta\omega_{\tau}^{\nu}} \right)
\overline{\vc{x}}^{\ast}(z)
&=& 
{\vc{x}(0)\vc{B}^{\ast}(\theta\omega_{\tau}^{\nu}) - \vc{x}(1)\vc{A}(0) 
\over \theta\omega_{\tau}^{\nu} - 1}
\nonumber
\\
&& \times
\lim_{z\to\theta\omega_{\tau}^{\nu}} 
\left(1 - {z\over \theta\omega_{\tau}^{\nu}} \right)
[\vc{I} - \vc{\Gamma}_A^{\ast}(z)]^{-1}.
\label{add-eqn-41}
\end{eqnarray}
We then obtain the following result.
\begin{lem}\label{lemma-limit-inverse-Gamma_A^{ast}(z)}
Suppose Assumptions \ref{assu-1}~(a)--(c), Assumption~\ref{assu-2}~(a)
and Assumption~\ref{assu-3} hold. Then for $\nu=0,1,\dots,\tau-1$,
\begin{equation}
\lim_{z \to \theta\omega_{\tau}^{\nu}}
\left(1 - {z \over \theta\omega_{\tau}^{\nu}} \right)
\left[ \vc{I} - \vc{\Gamma}_A^{\ast}(z) \right]^{-1} 
= {\vc{\Delta}_M(\omega_{\tau}^{\nu}) \vc{v}(\theta) 
\vc{\mu}(\theta)\vc{\Delta}_M(\omega_{\tau}^{\nu})^{-1}
\over (\rmd/\rmd y)\delta(\vc{A}^{\ast}(y))|_{y=\theta} - 1}.
\label{add-eqn-30}
\end{equation}
\end{lem}

\proof 
It follows from  
(\ref{eqn-inverse}), (\ref{eqn-det}) and (\ref{eq-adj-06}) that
\begin{eqnarray}
\lim_{z \to \theta\omega_{\tau}^{\nu}}
\left(1 - {z \over \theta \omega_{\tau}^{\nu}} \right)
\left[ \vc{I} - \vc{\Gamma}_A^{\ast}(z) \right]^{-1} 
&=&
{\lim_{z \to \theta \omega_{\tau}^{\nu}}}
{  1 - \dm{ z \over \theta \omega_{\tau}^{\nu}} 
 \over 1 - r_1(z) }
\cdot
\lim_{z \to \theta\omega_{\tau}^{\nu}}
{\adj( \vc{I} - \vc{\Gamma}_A^{\ast}(z) ) 
\over \displaystyle\prod_{m=2}^M (1 - r_m(z))}
\nonumber
\\
&=& { 1 \over \theta \omega_{\tau}^{\nu} }
{1 \over (\rmd / \rmd z) r_1(z) |_{z=\theta\omega_{\tau}^{\nu}} } 
\cdot
\vc{v}(\theta \omega_{\tau}^{\nu}) 
\vc{\mu}(\theta \omega_{\tau}^{\nu}),
\label{add-eqn-33}
\end{eqnarray}
where we use l'H\^{o}pital's rule in the second equality. Note that
(\ref{add-eqn-35}), (\ref{eqn-deff-delta}) and (\ref{def-r_1(z)})
yield
\[
\theta\omega_{\tau}^{\nu}
\left.
{\rmd \over \rmd z} r_1(z) \right|_{z=\theta\omega_{\tau}^{\nu}}
=  \theta \left. {\rmd \over \rmd z} r_1(z) \right|_{z=\theta}
= 
\left. {\rmd \over \rmd y} \delta(\vc{A}^{\ast}(y)) \right|_{y=\theta}-1,
~~~\nu=0,1,\dots,\tau-1.
\]
Thus (\ref{add-eqn-33}) leads to
\[
\lim_{z \to \theta\omega_{\tau}^{\nu}}
\left(1 - {z \over \theta\omega_{\tau}^{\nu}} \right)
\left[ \vc{I} - \vc{\Gamma}_A^{\ast}(z) \right]^{-1} 
= {\vc{v}(\theta \omega_{\tau}^{\nu}) 
\vc{\mu}(\theta \omega_{\tau}^{\nu})
\over (\rmd/\rmd y)\delta(\vc{A}^{\ast}(y))|_{y=\theta} - 1},
~~~\nu=0,1,\dots,\tau-1.
\]
Finally, substituting (\ref{add-eqn-36a}) and (\ref{add-eqn-36b}) into
the above equation, we obtain (\ref{add-eqn-30}).  \qed

\medskip

Applying Lemma \ref{lemma-limit-inverse-Gamma_A^{ast}(z)} to
(\ref{add-eqn-41}), we have
\begin{equation}
\lim_{z \to \theta\omega_{\tau}^{\nu}}
\left(1 - {z \over \theta\omega_{\tau}^{\nu}} \right)
\overline{\vc{x}}^{\ast}(z)
= c(\omega_{\tau}^{\nu})\cdot \vc{\mu}(\theta)
\vc{\Delta}_M(\omega_{\tau}^{\nu})^{-1},
\qquad \nu=0,1,\dots,\tau-1,
\label{add-eqn-31}
\end{equation}
where $c(\omega_{\tau}^{\nu})$ ($\nu=0,1,\dots,\tau-1$) is a scalar
such that
\begin{equation}
c(\omega_{\tau}^{\nu}) 
=
\dm{ \left[ 
\vc{x} (0) \vc{B}^{\ast}(\theta\omega_{\tau}^{\nu}) 
- \vc{x}(1) \vc{A}(0) 
\right]
\vc{\Delta}_M(\omega_{\tau}^{\nu})\vc{v}(\theta)
\over (\theta\omega_{\tau}^{\nu} - 1)
\{(\rmd/\rmd y)\delta(\vc{A}^{\ast}(y))|_{y=\theta} - 1\}
}.
\label{eqn-c(omega_{tau}^{nu})}
\end{equation}
It follows from (\ref{add-eqn-31}) that $c(\omega_{\tau}^{\nu}) \neq
0$ if and only if $z=\theta\omega_{\tau}^{\nu}$ is a simple pole of
$[\overline{\vc{x}}^{\ast}(z)]_j$ for any $j \in \M$.

\begin{lem}\label{lemma-c(1)}
Suppose Assumptions \ref{assu-1}--\ref{assu-3} hold and $\theta <
r_B$. Then $c(\omega_{\tau}^{0}) = c(1) > 0$ and therefore $z=\theta$
is a simple pole of $[\overline{\vc{x}}^{\ast}(z)]_j$ for any $j \in
\M$.
\end{lem}

\proof From (\ref{eqn-c(omega_{tau}^{nu})}), we have
\[
c(1)= {1 \over (\theta - 1)\{(\rmd/\rmd y)
\delta(\vc{A}^{\ast}(y))|_{y=\theta} - 1
\}}
\cdot
[\vc{x} (0) \vc{B}^{\ast}(\theta) 
- \vc{x}(1) \vc{A}(0) ]
\vc{v}(\theta).
\]
Thus according to (\ref{add-ineqn-01}) and
$\theta - 1 > 0$, it suffices to show that
\begin{equation}
[\vc{x} (0) \vc{B}^{\ast}(\theta) 
- \vc{x}(1) \vc{A}(0) ]
\vc{v}(\theta)
>  0.
\label{add-ineqn-02}
\end{equation}
From (\ref{eqn-x(k)-G}) and $\vc{G} = (\vc{I} -
\vc{U}(0))^{-1}\vc{A}(0)$, we have $\vc{x}(1)\vc{A}(0) =
\vc{x}(0)\vc{U}_0(1)\vc{G}$ and therefore
\begin{equation}
[\vc{x} (0) \vc{B}^{\ast}(\theta) 
- \vc{x}(1) \vc{A}(0) ]\vc{v}(\theta)
=
\vc{x} (0) \left(
\vc{B}^{\ast}(\theta) 
- \vc{U}_0(1)\vc{G}
\right)\vc{v}(\theta).
\label{add-eqn-46}
\end{equation}
Letting $z=\theta$ in (\ref{eqn-R_0}) and substituting it into
the right hand side of (\ref{add-eqn-46}), we have
\[
[\vc{x} (0) \vc{B}^{\ast}(\theta) 
- \vc{x}(1) \vc{A}(0) ]\vc{v}(\theta)
=
\vc{x}(0)\vc{R}_0^{\ast}(\theta)\vc{s}(\theta),
\]
where $\vc{s}(\theta) = (\vc{I}-\vc{U}(0))
(\vc{I}-\vc{G}/\theta)\vc{v}(\theta) \ge \vc{0}, \neq \vc{0}$ (see
Proposition~\ref{lemma5-Li05b}).  Note here that
$[\vc{R}_0^{\ast}(\theta)]_{i,j} > 0$ if and only if
$[\vc{R}_0^{\ast}(1)] > 0$. Note also that
$[\vc{R}_0^{\ast}(1)]_{i,j}$ represents the expected number of visits
to phase $j$ during the first passage time from state $(0,i)$ to level
zero. Thus $\vc{R}_0^{\ast}(\theta) \ge \vc{O}$ has no zero column due
to Assumption~\ref{assu-1}~(a). Consequently,
$\vc{x}(0)\vc{R}_0^{\ast}(\theta)\vc{s}(\theta) > 0$ because
$\vc{x}(0) > \vc{0}$. \qed

\begin{rem}\label{rem-Falk}
In Theorem~3.5 in \cite{Falk94}, (\ref{add-ineqn-02}) is assumed in
order to prove $c(1) > 0$.
\end{rem}

Lemma \ref{lemma-c(1)} ensures that (\ref{add-eqn-40}) holds. Thus
substituting (\ref{add-eqn-31}) into (\ref{add-eqn-40}) yields
\begin{equation}
\overline{\vc{x}}(k)
= \theta^{-k}\sum_{\nu=0}^{\tau-1} 
{1 \over (\omega_{\tau}^{\nu})^k}
c(\omega_{\tau}^{\nu})\cdot \vc{\mu}(\theta)
\vc{\Delta}_M(\omega_{\tau}^{\nu})^{-1}
+ O((\theta+\varepsilon_0)^{-k})\vc{e}^{\rm t}.
\label{asymp-overline{x}(k)}
\end{equation}
According to (\ref{add-eqn-57}), the dominant term of
(\ref{asymp-overline{x}(k)}) is positive for any $k=0,1,\dots$.
Letting $k=n\tau+l$ ($l=0,1,\dots,\tau-1$, $n=0,1,\dots$) in
(\ref{asymp-overline{x}(k)}), we readily obtain the following theorem.
\begin{thm}\label{main-thm-1}
Suppose Assumptions \ref{assu-1}--\ref{assu-3} hold and $\theta <
r_B$. Then
\[
\overline{\vc{x}}(n\tau+l) 
= \theta^{-n\tau-l} \vc{c}_l 
+ O((\theta+\varepsilon_0)^{-(n\tau+l)})\vc{e}^{\rm t},
\qquad l=0,1,\dots,\tau-1,
\]
where
\begin{equation}
\vc{c}_l
= 
\sum_{\nu=0}^{\tau-1}{1 \over (\omega_{\tau}^{\nu})^{l}}
\dm{ \left[ 
\vc{x} (0) \vc{B}^{\ast}(\theta\omega_{\tau}^{\nu}) 
- \vc{x}(1) \vc{A}(0) 
\right]
\vc{\Delta}_M(\omega_{\tau}^{\nu})\vc{v}(\theta)
\over (\theta\omega_{\tau}^{\nu} - 1)
\{(\rmd/\rmd y)\delta(\vc{A}^{\ast}(y))|_{y=\theta} - 1\}
}
\cdot 
\vc{\mu}(\theta)
\vc{\Delta}_M(\omega_{\tau}^{\nu})^{-1} > \vc{0}.
\label{def-c(l)}
\end{equation}
\end{thm}

\begin{rem}
Theorem \ref{main-thm-1} and Propositions \ref{prop-takine} and
\ref{prop-h-tau} imply that $\vc{c}_l$ ($l=0,1,\dots,\tau-1$) in
(\ref{def-c(l)}) must be equal to $\vc{d}_l$ in
Proposition~\ref{prop-takine}.  However in the case of $\tau = h \ge
2$, it seems difficult to demonstrate $\vc{c}_l = \vc{d}_l$
($l=0,1,\dots,\tau-1$), because these two expressions are completely
different. On the other hand, in the case of $\tau=h=1$, we can
readily confirm $\vc{c}_0$ is equal to the constant vector on the
right hand side of (\ref{add-eqn-50}) (see
Proposition~\ref{prop-takine-2}).
\end{rem}

\begin{coro}
Suppose Assumptions \ref{assu-1}--\ref{assu-3} hold and $\theta <
r_B$. If $\vc{C}(0) = \vc{A}(0)$, then
$\vc{c}_l$ ($l=0,1,\dots,\tau-1$) is given by
\begin{equation}
\vc{c}_l
=
\sum_{\nu=0}^{\tau-1}{1 \over (\omega_{\tau}^{\nu})^{l}}
\dm{ \vc{x} (0)\left[ 
 \vc{B}^{\ast}(\theta\omega_{\tau}^{\nu}) + \vc{B}(0) - \vc{I}
\right]
\vc{\Delta}_M(\omega_{\tau}^{\nu})\vc{v}(\theta)
\over (\theta\omega_{\tau}^{\nu} - 1)
\{(\rmd/\rmd y)\delta(\vc{A}^{\ast}(y))|_{y=\theta} - 1\} }
\cdot 
\vc{\mu}(\theta)
\vc{\Delta}_M(\omega_{\tau}^{\nu})^{-1}.
\label{eqn-c_l-01}
\end{equation}
In addition, if $\vc{B}(k) = \vc{A}(k)$ for all
$k=0,1,\dots$, then (\ref{eqn-c_l-01}) is reduced to
\begin{equation}
\vc{c}_l
=
\sum_{\nu=0}^{\tau-1}{1 \over (\omega_{\tau}^{\nu})^{l}}
\dm{ (1-\rho)\vc{g}
\vc{\Delta}_M(\omega_{\tau}^{\nu})\vc{v}(\theta)
\over 
(\rmd/\rmd y)\delta(\vc{A}^{\ast}(y))|_{y=\theta} - 1}
\cdot 
\vc{\mu}(\theta)
\vc{\Delta}_M(\omega_{\tau}^{\nu})^{-1},
\qquad l=0,1,\dots,\tau-1.
\label{eqn-c_l-02}
\end{equation}
\end{coro}

\proof When $\vc{C}(0) = \vc{A}(0)$, $\vc{x}(1)\vc{A}(0) = \vc{x}(0) -
\vc{x}(0)\vc{B}(0)$. Substituting this into (\ref{def-c(l)}) yields
(\ref{eqn-c_l-01}). We now suppose that $\vc{C}(0) = \vc{A}(0)$ and
$\vc{B}(k) = \vc{A}(k)$ for all $k=0,1,\dots$. We then have
\begin{eqnarray}
\left[
\vc{B}^{\ast}(\theta\omega_{\tau}^{\nu}) + \vc{B}(0) - \vc{I}
\right]
\vc{\Delta}_M(\omega_{\tau}^{\nu})\vc{v}(\theta)
&=& \left[\vc{A}^{\ast}(\theta\omega_{\tau}^{\nu}) - \vc{I}
\right]
\vc{\Delta}_M(\omega_{\tau}^{\nu})\vc{v}(\theta)
\nonumber
\\
&=& \left[(\theta\omega_{\tau}^{\nu})
\vc{\Gamma}_A^{\ast}(\theta\omega_{\tau}^{\nu}) - \vc{I}
\right]
\vc{\Delta}_M(\omega_{\tau}^{\nu})\vc{v}(\theta)
\nonumber
\\
&=& (\theta\omega_{\tau}^{\nu} - 1)
\vc{\Delta}_M(\omega_{\tau}^{\nu})\vc{v}(\theta),
\label{add-eqn-44}
\end{eqnarray}
where the second equality follows from
(\ref{add-eqn-Gamma_A^{ast}(z)}) and the third one follows from
$\vc{v}(\theta\omega_{\tau}^{\nu})=\vc{\Delta}_M(\omega_{\tau}^{\nu})\vc{v}(\theta)$
and $\delta(\vc{\Gamma}_A(\theta\omega_{\tau}^{\nu}))=1$ (see
(\ref{add-eqn-36b}) and Proposition \ref{prop-delta-Gamma_A-2}).
Substituting (\ref{add-eqn-44}) and $\vc{x}(0)=(1-\rho)\vc{g}$ (see,
e.g., Takine~\cite{Taki00}) into (\ref{eqn-c_l-01}), we have
(\ref{eqn-c_l-02}). \qed

\medskip

We close this subsection by discussing the period of the geometric
asymptotics of $\{\overline{\vc{x}}(k)\}$. It should be noted that
Theorem \ref{main-thm-1} does not necessarily show that the period in
the geometric asymptotics of $\{\overline{\vc{x}}(k)\}$ is equal to
$\tau$. This is because $c(\omega_{\tau}^{\nu})$
($\nu=1,2,\dots,\tau-1$) may be equal to zero i.e.,
$z=\theta\omega_{\tau}^{\nu}$ ($\nu=1,2,\dots,\tau-1$) may not be a
pole of $[\overline{\vc{x}}^{\ast}(z)]_j$'s ($j\in\M$). Let $\calP_A$
denote the set of poles of $[\overline{\vc{x}}^{\ast}(z)]_j$'s on $\{z
\in \C; |z|=\theta\}$. Since $c(1) > 0$ (see Lemma \ref{lemma-c(1)}),
$\calP_A$ includes $\theta$ and thus
\[
\calP_A = \{\theta \omega_{\tau}^{\nu}; \nu \in \bbH \},
\]
where $\bbH = \{0\} \cup \{\nu \in \{1,2,\dots,\tau-1\};
c(\omega_{\tau}^{\nu}) \neq 0\}$. Let $\tau' = \tau/\gcd\{\nu \in
\{\tau\} \cup \bbH\, \backslash \{0\}\}$. Let $H$ denote the number of
elements in $\bbH$. It is easy to see that there exists $H$
nonnegative integers, $\nu_m$'s ($m=0,1,\dots,H-1$), such that
$0=\nu_0 < \nu_1 < \cdots < \nu_{H-1} \le \tau'-1$ and
\begin{equation}
\calP_A
=
\{\theta\omega_{\tau'}^{\nu_m}; m=0,1,\dots,H-1\}.
\label{defn-calP_A}
\end{equation}
Note that since $\overline{\vc{x}}(k) > \vc{0}$ ($k=0,1,\dots$), each
$[\overline{\vc{x}}^{\ast}(z)]_j$ ($j \in \M$) has pairs of complex
conjugate poles and therefore
\begin{equation}
\omega_{\tau'}^{\nu_m} \omega_{\tau'}^{\nu_{H-m}} = 1,
\qquad m=1,2,\dots,\lfloor (H-1)/2 \rfloor.
\label{eqn-conjugate}
\end{equation}
It follows from (\ref{asymp-overline{x}(k)}) and
$c(\omega_{\tau}^{\nu}) = 0$ for $\nu \not\in \bbH$ that
\[
\overline{\vc{x}}(k)
= \theta^{-k}\sum_{m=0}^{H-1} {1 \over (\omega_{\tau'}^{\nu_m})^k}
c(\omega_{\tau'}^{\nu_m})\cdot \vc{\mu}(\theta)
\vc{\Delta}_M(\omega_{\tau'}^{\nu_m})^{-1}
+ O((\theta+\varepsilon_0)^{-k})\vc{e}^{\rm t}.
\]
Letting $k=n\tau'+l$ ($l=0,1,\dots,\tau'-1$, $n=0,1,\dots$) in the
above equation yields
\[
\overline{\vc{x}}(n\tau'+l) 
= \theta^{-n\tau'-l} \vc{c}_l' 
+ O((\theta+\varepsilon_0)^{-(n\tau'+l)})\vc{e}^{\rm t},
\qquad l=0,1,\dots,\tau'-1,
\]
where 
\begin{eqnarray*}
\vc{c}_l'
= 
\sum_{m=0}^{H-1}
{1 \over (\omega_{\tau'}^{\nu_m})^{l}}
\dm{ \left[ 
\vc{x} (0) \vc{B}^{\ast}(\theta\omega_{\tau'}^{\nu_m}) 
- \vc{x}(1) \vc{A}(0) 
\right]
\vc{\Delta}_M(\omega_{\tau'}^{\nu_m})\vc{v}(\theta)
\over (\theta\omega_{\tau'}^{\nu_m} - 1)
\{ (\rmd/\rmd y)\delta(\vc{A}^{\ast}(y))|_{y=\theta} - 1\} }
\nonumber
\\ 
{}
\cdot 
\vc{\mu}(\theta)
\vc{\Delta}_M(\omega_{\tau'}^{\nu_m})^{-1} > \vc{0}.
\end{eqnarray*}
As a result, the period in 
the geometric asymptotics of $\{\overline{\vc{x}}(k)\}$ is equal to
$\tau'$.

\subsection{Case of $\theta > r_B$}\label{subsec-case-2}

For simplicity, we denote, by $C(\zeta,r)$, the circle $\{z\in\C;
|z-\zeta| = r\}$ in the complex plane, where $\zeta \in \C$ and $r \ge
0$.  In this subsection, we make the following assumption.
\begin{assumption}\label{assu-B-ast}
(a) $\vc{B}^{\ast}(z)$ is meromorphic in an open set containing the
  domain $\{z\in\C;|z| \le r_B\}$\footnote{In the published version,
    it is assumed that ``$\vc{B}^{\ast}(z)$ is meromorphic in the
    domain $\{z\in\C;|z| \le r_B\}$". However, by the definition of
    meromorphicness, the revised description is more appropriate. The
    same is true of Assumption~\ref{assu-append-1}.}, and (b) there
  exists some positive integer $m_B$ and some finite nonnegative
  matrix $\widetilde{\vc{B}}(r_B)$ such that
\[
\lim_{z \to  r_B} \left( 
1 - {z \over r_B } \right)^{m_B}
\vc{B}^{\ast}(z)  = \widetilde{\vc{B}}(r_B) \neq \vc{O}.
\]
\end{assumption}

\begin{rem}
From the definition, $\vc{B}^{\ast}(z)$ is holomorphic in the domain
$\{z\in\C;|z| < r_B\}$. Thus Assumption \ref{assu-B-ast} is an
additional condition for the behavior of $\vc{B}^{\ast}(z)$ on the
convergence radius. In fact, Assumption \ref{assu-B-ast} shows that
$\vc{B}^{\ast}(z)$ is holomorphic on $C(0,r_B)$ except for its poles.
\end{rem}

It follows from Lemma \ref{lemma-append-2} that under
Assumption~\ref{assu-B-ast}, any pole of $[\vc{B}^{\ast}(z)]_{i,j}$'s
on $C(0,r_B)$ is of order less than or equal to $m_B$. Thus we assume
the following.
\begin{assumption}\label{add-assu-B-ast}
There exist exactly $N$ ($N \in \N$) complex numbers $\zeta_n$'s
$(n=0,1,\dots,N-1)$ on $C(0,1)$ such that $0 = \arg(\zeta_0) <
\arg(\zeta_1) < \cdots < \arg(\zeta_{N-1}) < 2\pi$ and
\[
\lim_{z \to  r_B \zeta_n} \left( 
1 - {z \over r_B \zeta_n} \right)^{m_B}
\vc{B}^{\ast}(z)  = \widetilde{\vc{B}}(r_B \zeta_n),
\qquad n=0,1,\dots,N-1,
\]
where $\widetilde{\vc{B}}(r_B \zeta_n)$ is some finite non-zero
matrix.
\end{assumption}

\begin{rem}\label{remark-conj}
Since $\vc{B}(k) \ge \vc{O}$ for all $k=1,2,\dots$,
\[
[\vc{B}^{\ast}(\conj(z))]_{i,j}
=
\conj\left([\vc{B}^{\ast}(z)]_{i,j} \right),
\qquad i, j \in \M,
\]
where $\conj(z)$ ($z \in \C$) denotes the complex conjugate of
$z$. Thus if $z=r_B\zeta$ is a pole of $[\vc{B}^{\ast}(z)]_{i,j}$, so is
$z=r_B\conj(\zeta)$. This fact implies that
\begin{equation}
\zeta_n\zeta_{N-n} = 1,
\qquad n=1,2,\dots, \lfloor (N-1) /2\rfloor,
\label{zeta-conjugate}
\end{equation}
and therefore for any $i,j \in \M$,
\[ 
\left[\widetilde{\vc{B}}(r_B\zeta_{N-n}) \right]_{i,j}
=
\conj\left([\widetilde{\vc{B}}(r_B \zeta_n)]_{i,j} \right),
\qquad n=1,2,\dots, \lfloor (N-1) /2\rfloor.
\]
\end{rem}

Lemma \ref{prop-smallest-pole} implies that if $\theta > r_B$, then
$\vc{I} - \vc{\Gamma}_A^{\ast}(z)$ is nonsingular in the domain
$\{z\in\C;|z| \le r_B\}$ and therefore $\overline{\vc{x}}^{\ast}(z)$
in (\ref{add-eqn-x^{ast}(z)-02}) is holomorphic in the same domain. It
thus follows from Assumptions~\ref{assu-B-ast} and
\ref{add-assu-B-ast} that for $ n=0,1,\dots,N-1$,
\begin{equation}
\lim_{z \to r_B \zeta_n} 
\left( 1 - {z \over r_B \zeta_n} \right)^{m_B}
\overline{\vc{x}}^{\ast}(z)
= {\vc{x}(0)\widetilde{\vc{B}}(r_B\zeta_n) \over r_B\zeta_n - 1}
\left[ \vc{I} - \vc{\Gamma}_A^{\ast}(r_B\zeta_n)) \right]^{-1}.
\label{add-eqn-34}
\end{equation}
Note here that for each $n=0,1,\dots,N-1$, $z=r_B\zeta_n$ is an
$m_B$th order pole of $[\overline{\vc{x}}^{\ast}(z)]_j$ ($j \in \M$)
if and only if
\[
\left[ 
\vc{x}(0)\widetilde{\vc{B}}(r_B\zeta_n)
\left[ \vc{I} - \vc{\Gamma}_A^{\ast}(r_B\zeta_n) \right]^{-1}
\right]_j \neq 0.
\]
Note also that
\[
\vc{x}(0)\widetilde{\vc{B}}(r_B)
\left[ \vc{I} - \vc{\Gamma}_A^{\ast}(r_B) \right]^{-1} > \vc{0},
\]
which follows from $\vc{x}(0) > \vc{0}$, $\widetilde{\vc{B}}(r_B)
\ge\vc{O},\neq \vc{O}$ and $[\vc{I} - \vc{\Gamma}_A^{\ast}(r_B)]^{-1}
> \vc{O}$. Therefore we obtain the following lemma.
\begin{lem}\label{lemma-pole-r_B}
If Assumptions \ref{assu-1}--\ref{assu-3} and \ref{assu-B-ast} hold
and $\theta > r_B$, then $z=r_B$ is an $m_B$th order pole of
$[\overline{\vc{x}}^{\ast}(z)]_j$ for any $j \in \M$.
\end{lem}
\begin{rem}\label{add-rem-pole-r_B}
Suppose all the conditions of Lemma~\ref{lemma-pole-r_B} hold except
Assumption \ref{assu-3}. Then in this case, we can readily show that
if $r_A > r_B$ (instead of $\theta>r_B$),
$\vc{x}(0)\widetilde{\vc{B}}(r_B) \left[ \vc{I} -
  \vc{\Gamma}_A^{\ast}(r_B) \right]^{-1} > \vc{0}$ and thus $z=r_B$ is
an $m_B$th order pole of $[\overline{\vc{x}}^{\ast}(z)]_j$ for any $j
\in \M$.
\end{rem}

From (\ref{add-eqn-34}) and Theorem~\ref{append-thm-asymp}, we readily
obtain the following result.
\begin{thm}
Suppose Assumptions \ref{assu-1}--\ref{assu-3}, \ref{assu-B-ast} and
\ref{add-assu-B-ast} hold and $\theta > r_B$. We then have
\begin{equation}
\overline{\vc{x}}(k) 
= {k^{m_B - 1} \over (m_B - 1)!}
{1 \over r_B^k} \vc{\xi}(k)
+ \left\{
\begin{array}{ll}
O((r_B + \varepsilon_0)^{-k})\vc{e}^{\rm t}, & m_B = 1,
\\
\rule{0mm}{5mm}
O(k^{m_B - 2} r_B^{-k})\vc{e}^{\rm t}, & m_B \ge 2,
\end{array}
\right. 
\label{add-eqn-55}
\end{equation}
where
\[
\vc{\xi}(k) 
= \sum_{n = 0 }^{N-1} 
 {1 \over \zeta_n^k}
{  \vc{x} (0)  \widetilde{\vc{B}}(r_B \zeta_n) \over r_B \zeta_n - 1}
\left[\vc{I} - \vc{\Gamma}_A^{\ast}(r_B\zeta_n) \right]^{-1}.
\]
Further $\limsup_{k\to\infty}\vc{\xi}(k)
> \vc{0}$ and $\vc{\xi}(k) \ge \vc{0}$ for all $k=0,1,\dots$.
\end{thm}

\begin{rem}
According to Theorem~\ref{append-thm-asymp}~(d), if $(\arg
\zeta_n)/\pi$ is rational for any $n=0,1,\dots,N-1$, then $\vc{\xi}(k)
> \vc{0}$ for any $k=0,1,\dots$.
\end{rem}

\begin{coro}
Suppose Assumptions~\ref{assu-1}--\ref{assu-3} and \ref{assu-B-ast}
hold and $\theta > r_B$. If $N = 1$, then
\begin{equation}
\overline{\vc{x}}(k) 
= {\vc{x}(0)\overline{\vc{B}}(k) \over r_B-1}
[\vc{I} - \vc{\Gamma}_A^{\ast}(r_B)]^{-1}
+ \left\{
\begin{array}{ll}
O((r_B + \varepsilon_0)^{-k})\vc{e}^{\rm t}, & m_B = 1,
\\
\rule{0mm}{5mm}
O(k^{m_B - 2} r_B^{-k})\vc{e}^{\rm t}, & m_B \ge 2.
\end{array}
\right.
\label{add-eqn-56}
\end{equation}
The dominant term on the right hand side of (\ref{add-eqn-56}) is a
positive vector for any $k=0,1,\dots$.
\end{coro}

We now mention the case where Assumption \ref{assu-3} does not hold,
i.e., there does not exist $\theta$ such that $1 < \theta < r_A$ and
$\theta = \delta(\vc{A}^{\ast}(\theta))$. In this case, $\det(\vc{I} -
\vc{\Gamma}_A^{\ast}(z)) = 0$ has no root in the domain $\{z\in\C;
1<|z|<r_A\}$, and thus if Assumption \ref{assu-B-ast} holds and $r_A >
r_B$, then $z=r_B$ is a dominant pole with order $m_B$ of
$[\overline{\vc{x}}^{\ast}(z)]_j$ for any $j \in \M$ (see Remark
\ref{add-rem-pole-r_B}). Therefore we have the following result.
\begin{thm}
Suppose there does not exist $\theta$ such that $1 < \theta < r_A$ and
$\theta = \delta(\vc{A}^{\ast}(\theta))$. Further suppose Assumptions
\ref{assu-1}, \ref{assu-2}, \ref{assu-B-ast} and \ref{add-assu-B-ast}
are satisfied and $r_A > r_B$. Then (\ref{add-eqn-55}) holds, and
(\ref{add-eqn-56}) does if $N=1$.
\end{thm}

\subsection{Case of $\theta = r_B$}\label{subsec-case-3}

This subsection considers the case of $\theta = r_B$ under Assumptions
\ref{assu-1}--\ref{assu-3}, \ref{assu-B-ast} and~\ref{add-assu-B-ast}.
Lemmas \ref{lemma-c(1)} and \ref{lemma-pole-r_B} show that $z=\theta$
(resp.~$r_B$) is a simple pole (resp.~an $m_B$th order pole) of each
$[\overline{\vc{x}}^{\ast}(z)]_j$ ($j \in \M$). Thus if $\theta =
r_B$, $z=\theta~(=r_B)$ is the $(m_B+1)$st order pole of each
$[\overline{\vc{x}}^{\ast}(z)]_j$ and its dominant poles are included
in $\calP \triangleq \calP_A \cap \calP_B$, where $\calP_A$ is given
in (\ref{defn-calP_A}) and $\calP_B =
\{\theta\zeta_n;n=0,1,\dots,N-1\}$.  Let $L$ denote the number of
elements in $\calP$. Let $\eta_m$ ($m=0,1,\dots,L-1$)'s denote $L$
nonnegative integers such that $0 = \eta_0 < \eta_1 <\cdots <
\eta_{L-1} \le \tau'-1$ and
\[
\calP = \{\theta\omega_{\tau'}^{\eta_m};m=0,1,\dots,L-1\}.
\]
Let $\widehat{\tau}$ denote
\[
\widehat{\tau} = \tau' / \gcd\{\eta_1,\eta_2,\dots,\eta_{L-1},\tau'\}.
\]
For simplicity, let $\widehat{\omega}_m =
(\omega_{\widehat{\tau}})^{\widehat{\eta}_m}$ ($m=0,1,\dots,L-1$), where
$\widehat{\eta}_m = (\widehat{\tau}/\tau')\eta_m$. It then follows
that $\calP=\{\theta\widehat{\omega}_m;m=0,1,\dots,L-1\}$. Note here
that $\widehat{\omega}_m\widehat{\omega}_{L-m} = 1$ for
$m=1,2,\dots,\lfloor(L-1)/2\rfloor$ due to (\ref{eqn-conjugate}) and
(\ref{zeta-conjugate}).

\begin{thm}\label{main-thm-3}
If Assumptions \ref{assu-1}--\ref{assu-3}, \ref{assu-B-ast} and
\ref{add-assu-B-ast} hold and $\theta = r_B$, then
\begin{equation}
\overline{\vc{x}} (n\widehat{\tau}+l) 
= { (n\widehat{\tau}+l)^{m_B} \over m_B!}
{1 \over \theta^{n\widehat{\tau}+l} } \widehat{\vc{c}}_l
+ O((n\widehat{\tau}+l)^{m_B-1} \theta^{-(n\widehat{\tau}+l)})\vc{e}^{\rm t},
\label{eq-rB-03}
\end{equation}
where
\[
\widehat{\vc{c}}_l
= \sum_{m=0}^{L-1} 
{1 \over (\widehat{\omega}_m)^l}
{ \vc{x}(0) 
\widetilde{\vc{B}}(\theta \widehat{\omega}_m) 
\vc{\Delta}_M(\widehat{\omega}_m)
\vc{v}(\theta)
\over (\theta \widehat{\omega}_m - 1) 
\{ (\rmd/\rmd y)\delta(\vc{A}^{\ast}(y))|_{y=\theta} - 1\} }
\cdot \vc{\mu}(\theta)
\vc{\Delta}_M(\widehat{\omega}_m)^{-1} > \vc{0}.
\]
\end{thm}

\begin{rem}
Theorem \ref{main-thm-3} shows that the period in the geometric
asymptotics of $\{\overline{\vc{x}}(k)\}$ is divisor of
$\widehat{\tau}$. It seems difficult to say more about the period in
the general setting.
\end{rem}

\noindent
{\it Proof of Theorem \ref{main-thm-3}}~ Recall that $z=\theta$ is an
$(m_B+1)$st order pole of $[\overline{\vc{x}}^{\ast}(z)]_j$ for any $j
\in \M$ and that $z=\theta\widehat{\omega}_m$ ($m=1,2,\dots,L-1)$ can
be minimum-modulus poles of order $m_B+1$. It then follows
from Theorem \ref{append-thm-asymp} that
\begin{equation}
\overline{\vc{x}}(k) 
=
\dm{ k^{m_B} \over m_B!} {1 \over \theta^k}
\sum_{m=0}^{L-1} 
{1 \over (\widehat{\omega}_m)^k}
\lim_{z \to \theta \widehat{\omega}_m} 
\left( 1 - {z \over \theta \widehat{\omega}_m} \right)^{m_B+1}
\overline{\vc{x}}^{\ast}(z)
+ O(k^{m_B-1} \theta^{-k})\vc{e}^{\rm t},
\label{eq-rB-04}
\end{equation}
where the dominant term is positive for any $k=0,1,\dots$.
Applying Assumption \ref{add-assu-B-ast} and Lemma
\ref{lemma-limit-inverse-Gamma_A^{ast}(z)} to
(\ref{add-eqn-x^{ast}(z)-02}) yield
\begin{eqnarray*}
\lefteqn{
\lim_{z \to \theta \widehat{\omega}_m} 
\left( 1 - {z \over \theta \widehat{\omega}_m} \right)^{m_B+1}
\overline{\vc{x}}^{\ast} (z) 
}
\qquad &&
\nonumber
\\
&=&
\lim_{z \to \theta \widehat{\omega}_m} 
\left( 1 - {z \over \theta \widehat{\omega}_m} 
\right)^{m_B}
{ \vc{x}(0) \vc{B}^{\ast}(z) - \vc{x}(1)\vc{A}(0)\over z - 1} 
\nonumber
\\
&& {} \cdot 
\lim_{z \to \theta \widehat{\omega}_m} 
\left( 1 - {z \over \theta \widehat{\omega}_m} 
\right)
[\vc{I} - \vc{\Gamma}_A^{\ast}(z)]^{-1}
\\
&=&
{ \vc{x}(0) 
\widetilde{\vc{B}}(\theta \widehat{\omega}_m) 
\over \theta \widehat{\omega}_m - 1}
{\vc{\Delta}_M(\widehat{\omega}_m)
\vc{v}(\theta) \cdot \vc{\mu}(\theta)
\vc{\Delta}_M(\widehat{\omega}_m)^{-1}
\over (\rmd/\rmd y)\delta(\vc{A}^{\ast}(y))|_{y=\theta} - 1},
\end{eqnarray*}
from which and (\ref{eq-rB-04}) we obtain
\begin{eqnarray*}
\overline{\vc{x}}(k) 
&=&
\dm{ k^{m_B} \over m_B!} {1 \over \theta^k}
\sum_{m=0}^{L-1} 
{1 \over (\widehat{\omega}_m)^k}
{ \vc{x}(0) 
\widetilde{\vc{B}}(\theta \widehat{\omega}_m) 
\over \theta \widehat{\omega}_m - 1}
{\vc{\Delta}_M(\widehat{\omega}_m)
\vc{v}(\theta) \cdot \vc{\mu}(\theta)
\vc{\Delta}_M(\widehat{\omega}_m)^{-1}
\over (\rmd/\rmd y)\delta(\vc{A}^{\ast}(y))|_{y=\theta} - 1}
\nonumber
\\
&& {} + O(k^{m_B-1} \theta^{-k})\vc{e}^{\rm t}.
\end{eqnarray*}
As a result, we obtain (\ref{eq-rB-03}) by letting
$k=n\widehat{\tau}+l$ ($n=0,1,\dots$, $l=0,1,\dots,\widehat{\tau}-1$)
in the above equation and using
$(\widehat{\omega}_m)^{n\widehat{\tau}+l} = (\widehat{\omega}_m)^l$.
\qed

\subsection{Remarks}

One of the referees informed the authors that a parallel research by
Dr.\ Tai \cite{Tai09} was open to the public after the submission of
this paper. The research is on the light-tailed asymptotics of the
stationary probability vectors $\{\vc{x}(k)\}$ of the Markov chain of
GI/G/1 type. Tai derives the decay rate of $\{\vc{x}(k)\}$ in a weak
sense, i.e., $-\log\lim_{k\to\infty}\left([\vc{x}(k)]_j\right)^{1/k}$,
and also presents several conditions under which $\{\vc{x}(k)\}$ is
asymptotically geometric, or light-tailed but not exactly geometric,
assuming the aperiodicity of the MAdP driven by the transition block
matrices in the non-boundary levels. As with this paper, Tai's
research includes the case where the jumps from the boundary level
have the dominant impact on the decay of the stationary tail
probability vectors.

\appendix

\section{Tail Asymptotics of Nonnegative Sequences}

Let $\{x_k;k=0,1,\dots\}$ denote a sequence of nonnegative numbers, an
infinite number of which are positive. Let $\sigma$ denote
\[
\sigma = \sup\left\{|z|; \sum_{k=0}^{\infty} x_k z^k < \infty, z \in \C\right\},
\] 
which is called the convergence radius of the power series.
Let $f(z)$ denote the generating function of $\{x_k;k=0,1,\dots\}$. We then have
\begin{equation}
f(z) = \sum_{k=0}^{\infty} x_k z^k,
\qquad |z| < \sigma.
\label{eqn-f(z)-MC}
\end{equation}
Further by definition, $f(z)$ is holomorphic inside the convergence radius.

In what follows, we make the following assumption.
\begin{assumption}\label{assu-append-1}
$f(z)$ is meromorphic in an open set containing the domain
  $\{z\in\C;|z| \le \sigma\}$\footnote{In the published version, it is
    assumed that ``$f(z)$ is meromorphic in the domain $\{z\in\C;|z|
    \le \sigma\}$".}, and the point $z=\sigma$ is an $\breve{m}$th
  pole of $f(z)$, where $\breve{m}$ is some finite positive integer.
\end{assumption}

\begin{lem}\label{lemma-append-2}
Under Assumption \ref{assu-append-1}, any pole of $f(z)$ on
$C(0,\sigma)$ is of order less than or equal to~$\breve{m}$.
\end{lem}

\proof 
We define $g(z)$ as
\[
g(z) = f(z) \left(1 - {z \over \sigma}\right)^{\breve{m}}.
\]
From (\ref{eqn-f(z)-MC}), we have for any $\varepsilon > 0$,
\begin{equation}
g(\sigma-\varepsilon) = f(\sigma-\varepsilon)\left({\varepsilon \over
\sigma}\right)^{\breve{m}} = \sum_{k=0}^{\infty} x_k (\sigma-\varepsilon)^k
\left({\varepsilon \over \sigma}\right)^{\breve{m}}.
\label{eqn-02}
\end{equation}
It thus follows from (\ref{eqn-f(z)-MC}) and (\ref{eqn-02}) that for
any $\omega_{\ast} \in \C$ such that $|\omega_{\ast}| = 1$ and
$\omega_{\ast} \neq 1$,
\begin{eqnarray}
\liminf_{\scriptstyle z=(\sigma-\varepsilon)\omega_{\ast} 
\atop \scriptstyle\varepsilon\downarrow 0} \left|
f(z) \left(1 - {z \over \sigma\omega_{\ast}}\right)^{\breve{m}}
\right|
&=& \liminf_{\varepsilon\downarrow 0}
\left|
\sum_{k=0}^{\infty} x_k (\sigma-\varepsilon)^k (\omega_{\ast})^k
\left({\varepsilon \over \sigma}\right)^{\breve{m}}
\right|
\nonumber
\\
&\le& \limsup_{\varepsilon\downarrow 0}
\sum_{k=0}^{\infty} x_k (\sigma-\varepsilon)^k 
\left({\varepsilon \over \sigma}\right)^{\breve{m}}
\nonumber
\\
&=& \limsup_{\varepsilon\downarrow 0}g(\sigma-\varepsilon)
= g(\sigma) < \infty,
\label{eqn-03}
\end{eqnarray}
where the last inequality holds because $g(z)$ is holomorphic in some
neighborhood of $z=\sigma$.  Let $\breve{m}_{\ast}$ denote
\[
\breve{m}_{\ast}
= \inf\left\{m \in \N \cup \{0\}; \lim_{z\to\sigma\omega_{\ast}}\left|
f(z) \left(1 - {z \over \sigma\omega_{\ast}}\right)^m
\right| < \infty
 \right\},
\]
where $f(z)(1-z/(\sigma\omega_{\ast}))^m$ is meromorphic in the domain
$\{z\in\C;|z| \le \sigma\}$ for $m = 0,1,\dots$. Thus, if
$\breve{m}_{\ast} > \breve{m}$, we have
\[
\liminf_{\scriptstyle z=(\sigma-\varepsilon)\omega_{\ast} 
\atop \scriptstyle\varepsilon\downarrow 0} 
\left|
f(z) \left(1 - {z \over \sigma\omega_{\ast}}\right)^{\breve{m}}
\right|
\ge  
\liminf_{z\to\sigma\omega_{\ast}}\left|
f(z)\left(1 - {z \over \sigma\omega_{\ast}}\right)^{\breve{m}}
\right| = \infty,
\]
which contradicts (\ref{eqn-03}). As a result, $\breve{m}_{\ast} \le
\breve{m}$, which implies that this lemma is true. \qed

\medskip

According to Lemma~\ref{lemma-append-2}, we introduce the following
definition.
\begin{defn}\label{defn-dominant-pole}
Under Assumption~\ref{assu-append-1}, {\it a dominant pole} of $f(z)$
is a pole that is located on its convergence circle $C(0,\sigma)$ and
is of the same order as that of pole $z=\sigma$. Thus the order of any
dominant pole of $f(z)$ is equal to $\breve{m}$.

\end{defn}

We make the following assumption, in addition to
Assumption~\ref{assu-append-1}.
\begin{assumption}\label{assu-append-2}
There exist exactly $P$ ($P \ge 1$) dominant poles, $\sigma_j$'s
($j=0,1,\dots,P-1$), of $f(z)$, where $\sigma_0=\sigma$ and $0 = \arg
\sigma_0 < \arg \sigma_1 < \cdots < \arg \sigma_{P-1} < 2\pi$.
\end{assumption}

\begin{rem}\label{appendix-rem}
Since $f(z)$ is the generating function of the nonnegative sequence
$\{x_k\}$, the set $\{\sigma_j;j=0,1,\dots,P-1\}$ consists of one or
two real numbers and $\lfloor (P-1)/2 \rfloor$ pairs of conjugate
complex numbers. Therefore $\sigma_j\sigma_{P-j} = \sigma^2$ for
$j=1,2,\dots,\lfloor (P-1)/2 \rfloor$.
\end{rem}

\begin{thm}\label{append-thm-asymp}
Suppose Assumptions \ref{assu-append-1} and \ref{assu-append-2} hold,
and let $a_{1,k} = 1/(\sigma+\varepsilon_0)^k$ ($k=0,1,\dots$), where
$\varepsilon_0 > 0$ is a sufficiently small number; and for
$m=2,3,\dots$, $a_{m,k} = k^{m-2} / \sigma^k$ ($k=0,1,\dots$). Then
the following are true.
\begin{enumerate}
\item The sequence $\{x_k;k=0,1,\dots\}$ satisfies
\begin{eqnarray}
x_k 
&=& {k+\breve{m}-1 \choose \breve{m}-1}{1 \over \sigma^k } 
\xi_k + O(a_{\breve{m},k})
\nonumber
\\
&=& {k^{\breve{m}-1} \over (\breve{m}-1)!}{1 \over \sigma^k }
\xi_k + O(a_{\breve{m},k}),
\label{asymp-x_k}
\end{eqnarray}
where
\begin{eqnarray}
\xi_k
&=& \sum_{j=0}^{P-1} \left({\sigma \over \sigma_j } \right)^k
 \lim_{z\to\sigma_j}
 \left(1 - {z \over \sigma_j} \right)^{\breve{m}} f(z).
\label{eq:coro-speed-1}
\end{eqnarray}
\item $\limsup_{k\to\infty}\xi_k > 0$.
\item $\xi_k \ge 0$ for all $k=0,1,\dots$.
\item In addition, if $\{x_k\}$ is eventually\footnote{In this revised
  version, we add ``eventually" in order to (slightly) strengthen the
  statement.} nonincreasing and $(\arg \sigma_j) / \pi $ is a rational
  number for any $j=0,1,\dots,P-1$, then $\xi_k > 0$ for all
  $k=0,1,\dots$.
\end{enumerate}
\end{thm}

\proof See Appendix \ref{proof-append-thm-asymp}.  \qed

\begin{rem}
A result similar to the statement (a) is given in Theorem~5.2.1
in \cite{Wilf94}.  Further when $P=1$, (\ref{eq:coro-speed-1}) is
reduced to eq.~(2) at p.~238 in \cite{Baio92}.
\end{rem}

\begin{rem}\label{rem-theorem-A.1}
Suppose the candidates for the dominant poles are $\sigma_j$'s
($j=0,1,\dots,P-1$) and at least one of them is indeed a dominant pole
(according to Assumption~\ref{assu-append-1}, $z=\sigma$ is a dominant
pole). For $\sigma_j$ not a dominant pole, we have
\[
\lim_{z\to\sigma_j}
 \left(1 - {z \over \sigma_j} \right)^{\breve{m}} f(z) = 0.
\]
Thus the statements (a)--(d) of Theorem~\ref{append-thm-asymp} still
hold, though the right hand side of (\ref{eq:coro-speed-1}) may
include some null terms.
\end{rem}

\section{Period of Markov Additive Processes}\label{append-MAdP}

This appendix summarizes fundamental results of the period of
MAdPs. In fact, most of the results described here are already implied
in \cite{Alsm94,Shur84}, though that is done in not an accessible
way. Further a MAdP related to the Markov chain of
M/G/1 type (and slightly more general one) are discussed in \cite{Gail96}.

We consider a MAdP $\{(\Gamma_n,J_n); n=0,1,\dots\}$, where the level
variable $\Gamma_n$ takes a value in $\Z = \{0, \pm 1, \pm 2, \dots\}$
and the phase variable $J_n$ takes a value in $\J \triangleq
\{1,2,\dots,J\}$. Let $\vc{\Gamma}(k)$ $(k\in\Z)$ denote a $J \times
J$ matrix whose $(i,j)$th $(i,j \in \J)$ element represents
\[
\Pr[\Gamma_{n+1}=k_0+k,J_{n+1}=j \mid \Gamma_n=k_0, J_n=i],
\]
for any fixed $k_0 \in \Z$. For simplicity, we denote the MAdP
$\{(\Gamma_n, J_n); n=0,1,\dots\}$ with kernel
$\{\vc{\Gamma}(k);k\in\Z\}$ by MAdP $\{\vc{\Gamma}(k);k\in\Z\}$. For
any two states $(k_1,j_1)$ and $(k_2,j_2)$ in $\Z \times \J$, we write
$(k_1,j_1) \rightarrow (k_2,j_2)$ when there exists a path from
$(k_1,j_1)$ to $(k_2,j_2)$ with some positive probability.
\begin{assumption}\label{assu-MAdP}
\hfill
\begin{enumerate}
\item $\vc{\Gamma} \triangleq \sum_{k\in \Z}\vc{\Gamma}(k)$ is
  irreducible.
\item For each $j \in \J$, there exists a nonzero integer $k_j$ such
  that $(0,j) \rightarrow (k_j,j)$.
\end{enumerate}
\end{assumption}

\medskip

Let $\K_j$
($j \in \J$) denote
\[
\K_j
=
\{ k \in \Z \backslash \{0\}; \; (0,j) \rightarrow (k,j)\},
\]
which is well-defined under Assumption \ref{assu-MAdP}.
\begin{lem}\label{add-lem-period}
Under Assumption \ref{assu-MAdP}, let $d_j = \gcd\{k\in\K_j\}$ for
$j\in\J$. Then $d_j$'s ($j \in \J$) are all
identical.
\end{lem}

\proof See Appendix~\ref{proof-lem-add-lem-period}. \qed

\begin{defn}\label{def-period-MAdP-01}
According to Lemma \ref{add-lem-period}, we write $d$ to represent
$d_j$'s and refer to the constant $d$ as the period of MAdP
$\{\vc{\Gamma}(k);k\in\Z\}$.
\end{defn}

We choose a state $i_0 \in \J$ and then define
$\J_0^{(i_0)}$ as
\[
\J_0^{(i_0)}
=
\{j \in \J; \; (0,i_0) \rightarrow (k,j), \; k \equiv 0~(\Mod~d)\}.
\]
We also define $\J_m^{(i_0)}$ ($m=1,2,\dots,d-1$) as
\[
\J_m^{(i_0)}
=
\{j \in \J ; \; (0,i_0) \rightarrow (k,j),
\; k \equiv m~(\Mod~d)\}.
\]
Since $\vc{\Gamma}$ is irreducible, each $j \in \J$ must belong to at
least one of $\{\J_m^{(i_0)};m=0,1,\dots,d-1\}$. Further for any $i
\in \J_m^{(i_0)}$,
\[
(0,i) \rightarrow (k,i_0)~\mbox{only if}~k \equiv -m~(\Mod~d),
\]
which implies that $\J_{m_1}^{(i_0)} \cap \J_{m_2}^{(i_0)} =
\emptyset$ for $m_1 \not\equiv m_2$ (if not, it would hold that
$(0,i_0) \rightarrow (k,i_0)$ for some $k \equiv m_1-m_2 \not\equiv
0~(\Mod~d)$). Thus $\J_0^{(i_0)}+\J_1^{(i_0)}+ \cdots
+\J_{d-1}^{(i_0)}=\J$ and there exists a function\footnote{corrected:
  ``an injective function" $\longrightarrow$ ``a function"} $q_0$ from
$\J$ to $\{0,1,\dots,d-1\}$ such that $j \in \J_{q_0(j)}^{(i_0)}$. It
follows from the definition of $\J_m^{(i_0)}$'s that
\[
[\vc{\Gamma}(k)]_{i,j} > 0
\mbox{ only if } k \equiv
q_0(j)-q_0(i)~(\Mod~d).
\]
As a result, we obtain the following result.
\begin{lem}\label{def-period-MAdP-02}
Under Assumption \ref{assu-MAdP}, the period $d$ is the largest
positive integer such that
\begin{equation}
[\vc{\Gamma}(k)]_{i,j} > 0
\mbox{ only if } k \equiv
q(j)-q(i)~(\Mod~d),
\label{eqn-add-30}
\end{equation}
where $q$ is some function\footnote{The published version states that
  function $q$ is injective. However, this is not true, in general.}
from $\J$ to $\{0,1,\dots,d-1\}$.  Further let $\J_m = \{j \in \J;
q(j)=m\}$ for $m=0,1,\dots,d-1$. Then $\J_m$'s ($m=0,1,\dots,d-1$) are
disjoint each other and $\J_{0}+\J_{1}+ \cdots +\J_{d-1}=\J$.
\end{lem}

In the rest of this section, we discuss the relationship between the
period $d$ of MAdP $\{\vc{\Gamma}(k);k\in\Z\}$ and the eigenvalues of
the generating function $\vc{\Gamma}^{\ast}(z)$ defined by $\sum_{k \in \Z}
z^k \vc{\Gamma}(k)$.  Let $\vc{\Delta}(z)$ denote a $J \times J$
diagonal matrix whose $j$th diagonal element is equal to
$z^{-q(j)}$. It then follows from (\ref{eqn-add-30}) that
\begin{equation}
\vc{\Gamma}^{\ast}(z)
= 
\vc{\Delta}(z)\vc{\Lambda}^{\ast}(z^d)
\vc{\Delta}(z)^{-1}
= 
\vc{\Delta}(z/|z|)\vc{\Lambda}^{\ast}(z^d)
\vc{\Delta}(z/|z|)^{-1},
\label{add-eqn-13}
\end{equation}
where $\vc{\Lambda}^{\ast}(z)$ denotes a $J \times J$ matrix whose
$(i,j)$th element is given by
\[
[\vc{\Lambda}^{\ast}(z)]_{i,j} = \sum_{n \in \Z}
z^n [\vc{\Gamma}(nd+q(j)-q(i))]_{i,j}.
\]
Let $\vc{\gamma}(z)$ and $\vc{g}(z)$ denote left- and
right-eigenvectors of $\vc{\Gamma}^{\ast}(z)$ corresponding to
eigenvalue $\delta(\vc{\Gamma}^{\ast}(z))$, normalized such that
\begin{equation}
\vc{\gamma}(z) \vc{\Delta}(z/|z|) \vc{e} = 1, 
\qquad
\vc{\gamma}(z)\vc{g}(z) = 1.
\label{cond-uniqueness}
\end{equation}
We then have the following lemma.
\begin{lem}\label{lem-delta-01}
Suppose Assumption \ref{assu-MAdP} holds, and let $I_{\gamma} = \{y >
0; \sum_{k \in \Z}y^k \vc{\Gamma}(k) < \infty\}$ and $\omega_x =
\exp(2\pi\iota/x)$ ($x\ge1$). Then the following hold for any $y \in
I_{\gamma}$ and $\nu=0,1,\dots,d-1$.
\begin{enumerate}
\item $\delta(\vc{\Gamma}^{\ast}(y\omega_d^{\nu})) =
  \delta(\vc{\Gamma}^{\ast}(y))$, both of which are simple
  eigenvalues.
\item $\vc{\gamma}(y\omega_d^{\nu}) =
  \vc{\gamma}(y)\vc{\Delta}(\omega_d^{\nu})^{-1}$ and
  $\vc{g}(y\omega_d^{\nu}) = \vc{\Delta}(\omega_d^{\nu})\vc{g}(y)$.
\end{enumerate}
\end{lem}

\proof It follows from (\ref{add-eqn-13}) that for $\nu =
0,1,\dots,d-1$,
\begin{eqnarray}
\vc{\Gamma}^{\ast}(y\omega_d^{\nu})
&=& 
\vc{\Delta}(y\omega_d^{\nu})\vc{\Lambda}^{\ast}(y^d)
\vc{\Delta}(y\omega_d^{\nu})^{-1},
\nonumber
\\
&=& \vc{\Delta}(\omega_d^{\nu})
[\vc{\Delta}(y)\vc{\Lambda}^{\ast}(y^d)
\vc{\Delta}(y)^{-1}]
\vc{\Delta}(\omega_d^{\nu})^{-1}
\nonumber
\\
&=& \vc{\Delta}(\omega_d^{\nu})
\vc{\Gamma}^{\ast}(y)
\vc{\Delta}(\omega_d^{\nu})^{-1},
\label{add-eqn-37}
\end{eqnarray}
which implies the statement (a) because $\vc{\Gamma}^{\ast}(y)$ is
nonnegative and irreducible. Next we prove the statement
(b). Pre-multiplying both sides of (\ref{add-eqn-37}) by
$\vc{\gamma}(y)\vc{\Delta}(\omega_d^{\nu})^{-1}$ and using
$\delta(\vc{\Gamma}^{\ast}(y\omega_d^{\nu})) =
\delta(\vc{\Gamma}^{\ast}(y))$, we have
\begin{eqnarray}
\left[ \vc{\gamma}(y)\vc{\Delta}(\omega_d^{\nu})^{-1} \right]
\vc{\Gamma}^{\ast}(y\omega_d^{\nu})
&=& \delta(\vc{\Gamma}^{\ast}(y)) 
\left[ \vc{\gamma}(y)\vc{\Delta}(\omega_d^{\nu})^{-1} \right]
\nonumber
\\
&=& \delta(\vc{\Gamma}^{\ast}(y\omega_d^{\nu})) 
\left[ \vc{\gamma}(y)\vc{\Delta}(\omega_d^{\nu})^{-1} \right].
\label{add-eqn-38}
\end{eqnarray}
Similarly we obtain
\begin{equation}
\vc{\Gamma}^{\ast}(y\omega_d^{\nu})
\left[\vc{\Delta}(\omega_d^{\nu})\vc{g}(y) \right]
= \delta(\vc{\Gamma}^{\ast}(y\omega_d^{\nu}))
\left[\vc{\Delta}(\omega_d^{\nu})\vc{g}(y) \right].
\label{add-eqn-39}
\end{equation}
It follows from (\ref{add-eqn-38}) and (\ref{add-eqn-39}) that there
exist some constants $\varphi_1$ and $\varphi_2$ such that
\[
\vc{\gamma}(y\omega_d^{\nu}) 
= \varphi_1 \vc{\gamma}(y) \vc{\Delta}(\omega_d^{\nu})^{-1},
\quad
\ \vc{g}(y\omega_d^{\nu}) 
= \varphi_2 \vc{\Delta}(\omega_d^{\nu}) \vc{g}(y).
\]
We can easily confirm that $\varphi_1 = \varphi_2 = 1$ satisfies the
normalizing condition (\ref{cond-uniqueness}). \qed

\begin{thm}\label{append-theorem}
Suppose Assumption~\ref{assu-MAdP} holds and
$\delta(\vc{\Gamma}^{\ast}(y)) = 1$ for some $y \in I_{\gamma}$, and let
$\omega$ denote a complex number such that $|\omega| = 1$. Then
$\delta(\vc{\Gamma}^{\ast}(y\omega)) = 1$ if and only if $\omega^d =
1$. Therefore
\begin{equation}
d = \max\{n \in \N; \delta(\vc{\Gamma}^{\ast}(y\omega_n)) = 1\}.
\label{add-eqn-52}
\end{equation}
Further if $\delta(\vc{\Gamma}^{\ast}(y\omega)) = 1$, the eigenvalue
is simple.
\end{thm}

\proof
Although Theorem~\ref{append-theorem} can be proved in a similar way
to Proposition~14 in \cite{Gail96}, the proof is given in
Appendix~\ref{proof-append-theorem} for completeness and the readers'
convenience. \qed

\begin{rem}
Theorem~\ref{append-theorem} provides a definition of the period of
MAdP $\{\vc{\Gamma}(k);k\in\Z\}$. In a very similar way, Shurenkov
\cite{Shur84} defined the period of MAdPs with proper kernels. In the
context of this paper, his definition is as follows:
\begin{equation}
d = \max\left\{n \in \N; \vc{\Gamma}^{\ast}(\omega_n)\vc{f} = \vc{f}
~\mbox{for some}~\vc{f} \in \C^J~\mbox{s.t.}~|[\vc{f}]_j| = 1~(j \in \J) \right\}.
\label{def-shurenkov}
\end{equation}
We can confirm that (\ref{def-shurenkov}) is equivalent to
Theorem~\ref{append-theorem} if $\vc{\Gamma}^{\ast}(1)$ is
stochastic. Shurenkov~\cite{Shur84} also implied that the statement of
Lemma \ref{def-period-MAdP-02} holds, based on which
Alsmeyer~\cite{Alsm94} defined the period of MAdPs.
\end{rem}

\section{Proofs}

\subsection{Proof of Proposition~\ref{prop-structure-G}}
\label{proof-prop-structure-G}

In order to prove this proposition by the reduction to absurdity, we
assume the negation of the statement, i.e., either (i) $\vc{G}$ is
strictly lower triangular, or (ii) $\vc{G}$ takes a form such that
\begin{equation}
\vc{G} 
=
\left(
\begin{array}{ccc}
\vc{G}_1     & \vc{O} &       \vc{O} \\
\vc{G}_{2,1} & \vc{G}_2     & \vc{O} \\
\vc{G}_{3,1} & \vc{G}_{3,2} & \vc{G}_3
\end{array}
\right),
\label{eq-G-form-add1}
\end{equation}
where $\vc{G}_i$ $(i= 1,2)$ is irreducible and $\vc{G}_2$ may be equal
to $\vc{G}_{\bullet}$ (in that case, $\vc{G}_{\bullet}$ is
irreducible). If case (i) is true, then $\vc{G}$ is a nilpotent matrix
and thus $\vc{G}^m = \vc{O}$ for some $m \in \N$, which is
inconsistent with Assumption~\ref{assu-1}~(a).

In what follows, we consider case (ii). For simplicity, we denote
$\{(k,j); j \in \M\}$ by $\bbL(k)$ $(k\in\N)$, and then partition
$\bbL(k)$ into subsets $\bbL_1 (k)$, $\bbL_2 (k)$ and $\bbL_3 (k)$
corresponding to $\vc{G}_1$, $\vc{G}_2$ and $\vc{G}_3$, respectively.
According to (\ref{eq-G-form-add1}) and the definition of $\vc{G}$,
for any $k \geq 2$ and $l$ $(1 \leq l < k)$ there exists no path from
$\bbL_1(k)$ to $\bbL_2(l)$ avoiding $\bbL(0) \triangleq \{(0,j); j \in
\M_0\}$.

We now fix $k \ge M+1\,(\ge 2)$.  Since the cardinality of $\M$ is
equal to $M$, it follows from Assumption~\ref{assu-1}~(b) that there
exists a path from $\bbL_1(k)$ to $\cup_{m=0}^{\infty}\bbL_2(m)$ of
length at most $M$. Such a path does not go through any state of
$\bbL(0)$ because of the skip-free-to-the-left property of
$\vc{T}$. Thus for some $l' \ge 1$, there exists a path from
$\bbL_1(k)$ to $\bbL_2(l')$ avoiding $\bbL(0)$. If $l' < k$, we
immediately have a contradiction. In fact, for $l' \ge k$, we also
have a contradiction because there exists a path from any state of
$\bbL_2(l')$ to any state of $\bbL_2(1)$ avoiding $\bbL(0)$, which
follows from the irreducibility of $\vc{G}_2$.  \qed

\subsection{Proof of Proposition \ref{prop-structure-R}}
\label{proof-prop-G-R}

We prove this proposition by reduction to absurdity. The negation of
the statement is that either (i) $\vc{R}$ is strictly upper triangular
or (ii) $\vc{R}$ takes a form such that
\begin{equation}
\vc{R} = 
\left(
\begin{array}{ccc}
\vc{R}_1 & \vc{R}_{1,2} & \vc{R}_{1,3}
\\
\vc{O}   & \vc{R}_2     & \vc{R}_{2,3}
\\
\vc{O}   &     \vc{O}   & \vc{R}_3
\end{array}
\right),
\label{eqn-partition-R}
\end{equation}
where $\vc{R}_i$ ($i=1,2$) is irreducible and $\vc{R}_2$ may be equal
to $\vc{R}_{\bullet}$. If case (i) is true, $\vc{R}$ is a nilpotent
matrix, which contradicts (\ref{eqn-delta(R^*(theta))}).

Next we consider case (ii). We partition $\bbL(k)$ into subsets
$\bbL_1(k)$, $\bbL_2(k)$ and $\bbL_3(k)$ corresponding to $\vc{R}_1$,
$\vc{R}_2$ and $\vc{R}_3$, respectively. Note here that\footnote{The
  equation of $[\vc{R}]_{i,j}$ is corrected.}
\[
[\vc{R}]_{i,j} 
= \E\left[\left. \sum_{k=1}^{\infty}
\sum_{n=1}^{a(k)-1} 1(X_n=k+1,~S_n = j) 
\,\right|\, X_0 = 1, S_0 = i\right],
\]
where $1(\chi)$ denotes the indicator function of event $\chi$. Thus
(\ref{eqn-partition-R}) implies that for any $l \ge 2$ there exists no
path from $\bbL_2(1)$ to $\bbL_1(l)$ avoiding $\bbL(1)$. On the other
hand, owing to the irreducibility of $\vc{R}_2$, there exists some
integer $k_{\ast} \ge M+2$ such that $\bbL_2(k_{\ast})$ is reachable
from $\bbL_2(1)$ avoiding $\bbL(1)$. Further for some $l_{\ast} \ge
2$, there exists a path from $\bbL_2(k_{\ast})$ to $\bbL_1(l_{\ast})$
avoiding $\bbL(1)$, because the cardinality of $\M$ is equal to $M$,
$\vc{A}$ is irreducible and $\vc{T}$ has the skip-free-to-the-left
property. Therefore we have a path from $\bbL_2(1)$ to
$\bbL_1(l_{\ast})$ via $\bbL_2(k_{\ast})$ ($l_{\ast},k_{\ast} \ge 2$)
avoiding $\bbL(1)$. This yields a contradiction. \qed

\subsection{Proof of Proposition \ref{prop-period}}
\label{proof-lem-period}

We suppose that there exists some $i \in \M$ such that $(0,i)
\rightarrow (k,i)$ only for $k=0$. Since
$\sum_{k\in\Z}\vc{\Gamma}_A(k)=\vc{A}$ is irreducible (see Assumption
\ref{assu-1}~(b)), for each $j \in \M$ there exists a unique pair
$(k_{i,j},k_{j,i})$ such that $k_{i,j} + k_{j,i} = 0$ and
\[
(0,i) \rightarrow (k_{i,j},j),
\qquad
(0,j) \rightarrow (k_{j,i},i).
\]
It thus follows that for any $k_0 \in \Z$,
\begin{equation}
\Pr[ \breve{X}_n \ge k_0 + K_{\rm min}~(\forall n=1,2,\dots) 
\mid \breve{X}_0 = k_0, \breve{S}_0 = i]
= 1,
\label{eqn-max-level}
\end{equation}
where $K_{\rm min} = \min_{j \in \M} k_{i,j}$.  We now fix $k_0$ to
be $k_0 = \max(1,1-K_{\rm min})$. Clearly $k_0 \ge 1$ and $k_0 +
K_{\rm min} = \max(1,1+K_{\rm min}) \ge 1$. Therefore
(\ref{eqn-max-level}) yields
\begin{equation}
\Pr[\breve{X}_n \ge 1~(\forall n=1,2,\dots) 
\mid \breve{X}_0=k_0, \breve{S}_0=i]=1.
\label{eqn-prob-equal}
\end{equation}
As a result, from (\ref{add-eqn-58}) and (\ref{eqn-prob-equal}), we
have
\[
\Pr[X_n \ge 1~(\forall n=1,2,\dots) \mid X_0=k_0, S_0=i]=1,
\]
which contradicts Assumption \ref{assu-1}~(a). \qed

\subsection{Proof of Lemma \ref{lem-adj}}
\label{proof-lem-adj}

It follows from the definition of $r_m(z)$'s ($m=1,2,\dots,M$) that
\begin{equation}
\adj(x\vc{I} - \vc{\Gamma}_A^{\ast}(z)) 
( x\vc{I} - \vc{\Gamma}_A^{\ast}(z))
=
\prod_{m=1}^{M} (x - r_m(z)) \vc{I}.
\label{eq-adj-02}
\end{equation}
Since $r_1(\theta\omega_{\tau}^{\nu}) =
\delta(\vc{\Gamma}^{\ast}_A(\theta\omega_{\tau}^{\nu})) = 1$
($\nu=0,1,\dots,\tau-1$), (\ref{eq-adj-02}) leads to
\begin{equation}
\adj ( \vc{I} - \vc{\Gamma}_A^{\ast}(\theta \omega_{\tau}^{\nu})) 
( \vc{I} - \vc{\Gamma}_A^{\ast}(\theta \omega_{\tau}^{\nu})) = \vc{O},
\qquad \nu=0,1,\dots,\tau-1.
\label{eq-adj-04}
\end{equation}
On the other hand, differentiating both sides of (\ref{eq-adj-02})
with respect to $x$ yields
\[
\left[
{\rmd \over \rmd x}\adj(x\vc{I} - \vc{\Gamma}_A^{\ast}(z)) 
\right]
(x\vc{I} - \vc{\Gamma}_A^{\ast}(z))
+ \adj(x\vc{I} - \vc{\Gamma}_A^{\ast}(z))
= \sum_{l=1}^M\prod_{m\in\M\backslash\{l\}} (x - r_m(z)) \vc{I}.
\]
Post-multiplying both sides of the above equation by
$\vc{v}(z)\vc{\mu}(z)$, we obtain
\begin{eqnarray}
&&\left[
{\rmd \over \rmd x}\adj(x\vc{I} - \vc{\Gamma}_A^{\ast}(z)) 
\right]
[x - r_1(z)]\vc{v}(z)\vc{\mu}(z)
+ \adj(x\vc{I} - \vc{\Gamma}_A^{\ast}(z))\vc{v}(z)\vc{\mu}(z)
\nonumber
\\
&& \qquad 
=  \sum_{l=1}^M\prod_{m\in\M\backslash\{l\}} (x - r_m(z))\vc{v}(z)\vc{\mu}(z).
\label{diff-eq-adj-02}
\end{eqnarray}
Since $r_1(\theta\omega_{\tau}^{\nu}) = 1$ and
$\prod_{m=2}^M(1-r_m(\theta\omega_{\tau}^{\nu})) \neq 0$, letting
$x=1$ and $z = \theta\omega_{\tau}^{\nu}$ in (\ref{diff-eq-adj-02})
yields
\[
\adj ( \vc{I} - \vc{\Gamma}_A^{\ast} (\theta \omega_{\tau}^{\nu})) 
\vc{v}(\theta \omega_{\tau}^{\nu}) \vc{\mu} (\theta \omega_{\tau}^{\nu})
=
\prod_{m=2}^{M} (1 - r_m(\theta \omega_{\tau}^{\nu})) 
\vc{v} (\theta \omega_{\tau}^{\nu}) \vc{\mu}(\theta \omega_{\tau}^{\nu}),\]
from which and (\ref{eq-adj-04}) it follows that
\begin{eqnarray}
&&
\adj ( \vc{I} - \vc{\Gamma}_A^{\ast}(\theta \omega_{\tau}^{\nu})) 
[ \vc{I} - \vc{\Gamma}_A^{\ast}(\theta \omega_{\tau}^{\nu}) 
+ \vc{v} (\theta \omega_{\tau}^{\nu})
\vc{\mu} (\theta \omega_{\tau}^{\nu})]
\nonumber
\\
&& {} \qquad 
=
\prod_{m=2}^{M} (1 - r_m(\theta \omega_{\tau}^{\nu})) 
\vc{v} (\theta \omega_{\tau}^{\nu}) \vc{\mu} (\theta \omega_{\tau}^{\nu}).
\label{eq-adj-05}
\end{eqnarray}
Note here that $\vc{I} - \vc{\Gamma}_A^{\ast}(\theta
\omega_{\tau}^{\nu}) + \vc{v} (\theta \omega_{\tau}^{\nu}) \vc{\mu}
(\theta \omega_{\tau}^{\nu})$ is non-singular and
\[
\vc{\mu} (\theta \omega_{\tau}^{\nu}) 
\left[ \vc{I} - \vc{\Gamma}_A^{\ast}(\theta \omega_{\tau}^{\nu}) 
+ \vc{v} (\theta \omega_{\tau}^{\nu}) 
\vc{\mu} (\theta \omega_{\tau}^{\nu})\right]^{-1} 
= \vc{\mu}(\theta \omega_{\tau}^{\nu}).
\]
Thus (\ref{eq-adj-05}) leads to (\ref{eq-adj-06}). \qed

\subsection{Proof of Theorem \ref{append-thm-asymp}}
\label{proof-append-thm-asymp}

\quad {\it Statement (a).}  It follows from
Assumption~\ref{assu-append-1} that there exists some $R > \sigma$
such that $f(z)$ is holomorphic in the domain $\{z\in\C; \sigma < |z|
\le R\}$.  We can choose $P$ positive numbers $r_j$'s
($j=0,1,\dots,P-1$) such that all the $C(\sigma_j,r_j)$'s are strictly
inside $C(0,R)$ and any two of them have no intersection. Let $\D$
denote
\[
\D = \{z; |z| < R,~|z - \sigma_j| > r_j,~j=0,1,\dots,P-1\}. 
\]
Clearly $f(z)$ is holomorphic in domain $\D \cup C(0,R) \cup
C(\sigma_0,r_0) \cup \cdots \cup C(\sigma_{P-1},r_{P-1})$. Thus by the
Cauchy integral formula, we have
\begin{equation}
f(z) =
{1 \over 2\pi \iota} \oint_{C(0,R)} 
{ f(\zeta ) \over \zeta - z} \rmd\zeta
-
{1 \over 2\pi \iota} \sum_{j=0}^{P-1}
\oint_{C(\sigma_{j},r_{j})} { f(\zeta ) \over \zeta - z} \rmd\zeta,
\qquad z \in \D,
\label{eq:f-co}
\end{equation}
where the integrals are taken counter-clockwise. 

We now consider the first
term in (\ref{eq:f-co}). For any $z \in \D$ and $\zeta \in
C(0,R)$, we have $|z / \zeta| < 1$ and therefore
\begin{equation}
{1 \over 2\pi \iota} \oint_{C(0,R)}
{ f(\zeta ) \over \zeta - z} \rmd\zeta 
= {1 \over 2\pi \iota} \oint_{C(0,R)}{f(\zeta) \over \zeta}
\sum_{n=0}^{\infty} { z^{n} \over \zeta^n} \rmd\zeta,
\qquad z \in \D.
\label{add-eqn-28}
\end{equation}
Since $f(\zeta)$ is holomorphic for $\zeta \in C(0,R)$, there exists
some $f_{\rm max} > 0$ such that
\begin{equation}
|f(\zeta)| \le f_{\rm max},
\qquad \zeta \in C(0,R).
\label{bound-f(zeta)}
\end{equation}
Thus for any fixed $z \in \D$,
\[
\left| {f(\zeta) \over \zeta}
\sum_{n=0}^{\infty} { z^{n} \over \zeta^n} \right|
\le {f_{\rm max} \over R} 
\sum_{n=0}^{\infty}
\left| { z  \over R} \right|^{n}
= {f_{\rm max}\over R} {1 \over 1 - \left|  \dm{z  \over R} \right|}
< \infty,
\qquad \zeta \in C(0,R),
\]
which shows that the order of summation and integration on the right
hand side of (\ref{add-eqn-28}) is interchangeable. As a result, it
follows from (\ref{add-eqn-28}) that
\begin{equation}
{1 \over 2\pi \iota} 
\oint_{C(0,R)} { f(\zeta ) \over \zeta - z} \rmd\zeta
=
\sum_{n=0}^{\infty} \left( 
{1 \over 2\pi \iota} \oint_{C(0,R)}
{ f(\zeta )  \over \zeta^{n+1}} \rmd\zeta
\right) z^{n}
= \sum_{n=0}^{\infty} c_n z^n,
\label{eq:f-1}
\end{equation}
where
\begin{equation}
c_{n} 
=
\displaystyle{{1 \over 2\pi \iota} \oint_{C(0,R)}}
{ f(\zeta )  \over \zeta ^{n+1}} \rmd\zeta,
\qquad
n=0,1,\dots.
\label{eqn-c_k}
\end{equation}

Next we consider the second term in (\ref{eq:f-co}). Since $|\zeta -
\sigma_{j}| / |z -\sigma_{j}| < 1$ for any $z \in \D$ and $\zeta
\in C(\sigma_{j},r_{j})$,
\[
{1 \over \zeta - z} 
=
{1 \over (\sigma_{j} - z) - (\sigma_{j} - \zeta) } 
=
{1 \over \sigma_{j} - z} 
\cdot 
{1 \over 1 - 
\displaystyle{ 
{\sigma_{j} - \zeta \over \sigma_{j} - z } 
}
}
=
{1 \over \sigma_{j} - z} \sum_{n=0}^{\infty} 
\left( {\sigma_{j} - \zeta \over \sigma_{j} - z } \right)^n.
\]
Thus we have
\[
{1 \over 2\pi \iota} 
\oint_{C(\sigma_{j},r_{j})}
{ f(\zeta ) \over \zeta - z} \rmd\zeta 
=
 {1 \over 2\pi \iota} 
\oint_{C(\sigma_{j},r_{j})}
\sum_{n=1}^{\infty} { f(\zeta ) (\sigma_{j} - \zeta)^{n-1} 
\over (\sigma_{j} - z)^n} \rmd\zeta,
\qquad z \in \D.
\]
In a way very similar to the right hand side of (\ref{add-eqn-28}), we
can confirm that the order of summation and integration in the above
equation is interchangeable, and then obtain
\begin{eqnarray}
{1 \over 2\pi \iota} 
\oint_{C(\sigma_{j},r_{j})} { f(\zeta ) \over \zeta - z} \rmd\zeta 
=
\sum_{n=1}^{\infty}(-1)^{n-1}
\left(
{1 \over 2\pi \iota} \oint_{C(\sigma_{j},r_{j})}
  f(\zeta )  (\zeta - \sigma_{j})^{n-1} \rmd\zeta
\right)
\nonumber
\\
\cdot {} { 1 \over (\sigma_{j} - z)^{n}},
\qquad z \in \D. \qquad 
\label{add-eqn-29}
\end{eqnarray}
Since $z=\sigma_{j}$ is an $\breve{m}$th order
pole,
\[
{1 \over 2\pi \iota} \oint_{C(\sigma_{j},r_{j})}
  f(\zeta )  (\zeta - \sigma_{j})^{n-1} \rmd\zeta = 0,
\qquad \mbox{for all}~n=\breve{m}+1,\breve{m}+2,\dots,
\]
from which and (\ref{add-eqn-29}) we have
\begin{eqnarray}
\lefteqn{
{1 \over 2\pi \iota} 
\oint_{C(\sigma_{j},r_{j})} { f(\zeta ) \over \zeta - z} \rmd\zeta 
}\qquad &&
\nonumber
\\
&=&
\sum_{n=1}^{\breve{m}}(-1)^{n-1}
\left(
{1 \over 2\pi \iota} \oint_{C(\sigma_{j},r_{j})}
  f(\zeta )  (\zeta - \sigma_{j})^{n-1} \rmd\zeta
\right){ 1 \over (\sigma_{j} - z)^{n}}
\nonumber
\\
&=& -
\sum_{n=1}^{\breve{m}} \sigma_{j}^{\breve{m}}
\left[
(-1)^{n+\breve{m}} 
\left(
{1 \over 2\pi \iota} \oint_{C(\sigma_{j},r_{j})}
 { f(\zeta ) (1 - \zeta/\sigma_{j})^{\breve{m}} 
\over (\zeta - \sigma_{j})^{\breve{m}-n+1}}\rmd\zeta
\right) 
\right]
{ 1 \over (\sigma_{j} - z)^{n}}
\nonumber
\\
&=& -
\sum_{n=1}^{\breve{m}}{\sigma_{j}^{\breve{m}} c_{j,n} \over (\sigma_{j} - z)^{n}},
\label{add-eqn-32}
\end{eqnarray}
where
\begin{equation}
c_{j,n}
=
(-1)^{n+\breve{m}}
\cdot
{1 \over 2\pi \iota} \oint_{C(\sigma_{j},r_{j})}
 { f(\zeta ) (1 - \zeta/\sigma_{j})^{\breve{m}} 
\over (\zeta - \sigma_{j})^{\breve{m}-n+1}} \rmd\zeta.
\label{def-c_{j,-n}}
\end{equation}
Substituting (\ref{eq:f-1}) and (\ref{add-eqn-32}) into
(\ref{eq:f-co}), we have
\[
f(z) 
= 
\sum_{n=0}^{\infty} c_n z^n
+
\sum_{j=0}^{P-1}
\sum_{n=1}^{\breve{m}} {\sigma_{j}^{\breve{m}} c_{j,n} \over ( \sigma_{j} - z)^{n}},
\qquad z \in \D,
\]
and therefore
\begin{equation}
\sum_{n=0}^{\infty}x_n z^n
=
\sum_{n=0}^{\infty} c_n z^n
+
\sum_{j=0}^{P-1}
\sum_{n=1}^{\breve{m}} { \sigma_{j}^{\breve{m}} c_{j,n} \over ( \sigma_{j} - z)^{n}},
\qquad z \in \D \cap \{ z \in \C; |z| < \sigma\}.
\label{eq:laurent}
\end{equation}
Differentiating both sides of (\ref{eq:laurent}) $k$ times with
respect to $z$, dividing them by $k!$ and letting $z=0$ yield
\begin{equation}
x_k
= 
c_k
+
\sum_{j=0}^{P-1}
\sum_{n=1}^{\breve{m}} \sigma_{j}^{\breve{m}-n} c_{j,n}
\dm{k+n-1 \choose n-1} 
{ 1 \over \sigma_{j}^k }.
\label{eqn-partial-fract-x_n}
\end{equation}
It follows from (\ref{bound-f(zeta)}) and (\ref{eqn-c_k}) that
\[
|c_k| 
\le {1 \over 2\pi} 
\oint_{C(0,R)} \left|{f(\zeta) \over \zeta^{k+1}} \right| \rmd\zeta 
\le {1 \over 2\pi} \oint_{C(0,R)} {f_{\rm max} \over R^{k+1}} \rmd\zeta
= {f_{\rm max} \over R^k},
\]
which leads to
\begin{equation}
\lim_{k\to\infty}
\left|{c_k \over \dm{1 \over \sigma_j^{k} }} \right|
= \lim_{k\to\infty}
|c_k|  \sigma^{k}
\le
\lim_{k\to\infty} f_{\rm max} 
\left( { \sigma \over R} \right)^{k}
= 0,
\qquad \mbox{for all}~j=0,1,\dots,P-1,
\label{add-eqn-51}
\end{equation}
where we use $|\sigma_j|=\sigma$ ($j=0,1,\dots,P-1$) and $0 < \sigma/R
< 1$. From (\ref{eqn-partial-fract-x_n}) and (\ref{add-eqn-51}), we
have
\begin{eqnarray}
x_k
&=& \dm{k+\breve{m}-1 \choose \breve{m}-1}{1 \over \sigma^k}
\sum_{j=0}^{P-1}   
\left({\sigma \over \sigma_{j}}\right)^k c_{j,\breve{m}}
+ O(a_{\breve{m},k})
\nonumber
\\
&=&
\dm{k^{\breve{m}-1} \over (\breve{m}-1)!} {1 \over \sigma^k}
\sum_{j=0}^{P-1} \left({\sigma \over \sigma_{j}}\right)^k c_{j,\breve{m}}
+ O(a_{\breve{m},k}).
\label{add-eqn-A00}
\end{eqnarray}
Note here that (\ref{def-c_{j,-n}}) yields
\begin{equation}
c_{j,\breve{m}}
=
{1 \over 2\pi \iota} \oint_{C(\sigma_{j},r_{j})}
 { f(\zeta ) (1 - \zeta/\sigma_{j})^{\breve{m}} 
\over \zeta - \sigma_{j}}\rmd\zeta
= 
\lim_{\zeta\to\sigma_j}
\left(1 - {\zeta \over \sigma_{j}} \right)^{\breve{m}} f(\zeta),
\label{add-eqn-A01}
\end{equation}
where we use the Cauchy integral formula in the last equality. As a
result, the statement (a) is true.

\medskip

\noindent
\quad {\it Statement (b). }From (\ref{asymp-x_k}) and the definition
of $\{a_{\breve{m},k}\}$, we have
\begin{equation}
x_k = \dm{k^{\breve{m}-1} \over (\breve{m}-1)!} {1 \over \sigma^k}\xi_k
+ o\left({k^{\breve{m}-1} \over \sigma^k}\right).
\label{add-eqn-A04}
\end{equation}
We now suppose $\limsup_{k\to\infty}\xi_k \le 0$.  Then
(\ref{add-eqn-A04}) yields
\[
\limsup_{k\to\infty}{x_k \over k^{{\breve{m}}-1}\sigma^{-k}}
= 0,
\]
which implies that for any $\varepsilon > 0$ there exists some
positive integer $K_{\varepsilon} \ge \breve{m}-1$ such that $x_k <
\varepsilon(k^{\breve{m}-1}/ \sigma^k)$ for all $k = K_{\varepsilon},
K_{\varepsilon}+1,\dots$. Thus we have
\begin{equation}
f(y) 
\le \sum_{k=0}^{K_{\varepsilon}-1} y^k x_k 
+ \varepsilon 
\sum_{k=K_{\varepsilon}}^{\infty}k^{\breve{m}-1}\left({y \over \sigma}\right)^k,
\qquad 0 \le y < \sigma.
\label{add-eqn-A05}
\end{equation}
Note that for $l=1,2,\dots$,
\begin{equation}
\sum_{k=l}^{\infty}k(k-1)\cdots(k-l+1)\left({y \over \sigma}\right)^k
= (-1)^{l+1} l! {\sigma y^l \over (y - \sigma)^{l+1}}.
\label{add-eqn-A06}
\end{equation}
Note also that there exists an $(\breve{m}-1)$-tuple
$(b_1,b_2,\dots,b_{\breve{m}-1})$ of real numbers such that
\begin{equation}
k^{\breve{m}-1} = \sum_{l=1}^{\breve{m}-1}b_l \cdot k(k-1) \cdots (k-l+1).
\label{add-eqn-A07}
\end{equation}
It follows from (\ref{add-eqn-A05}), (\ref{add-eqn-A06}) and
(\ref{add-eqn-A07}) that for any $\varepsilon > 0$,
\[
0 \le \limsup_{y\uparrow\sigma}\left(1 - {y \over \sigma} \right)^{\breve{m}}f(y) 
\le \varepsilon  b_{\breve{m}-1} (\breve{m}-1)! .
\]
Letting $\varepsilon \to 0$ in the above inequality, we have
$\lim_{y\uparrow\sigma}\left\{1 - (y / \sigma)
\right\}^{\breve{m}}f(y) = 0$, which is inconsistent with
Assumption~\ref{assu-append-1}.

\medskip
\noindent
\quad {\it Statement (c).} 
It follows from (\ref{add-eqn-A01}), Assumption~\ref{assu-append-2}
and Remark~\ref{appendix-rem} that $c_{0,\breve{m}}$ is a real number
and ($c_{j,\breve{m}},c_{P-j,\breve{m}}$) ($j=1,2,\dots,\lfloor (P-1)/2
\rfloor$) is a pair of complex conjugates, and thus $\xi_k$ is a real
number such that
\begin{equation}
\xi_k = y_0 
+ \sum_{j=1}^{\lfloor (P-1)/2\rfloor} y_j \cos(2\pi k\alpha_j),
\qquad k=0,1,\dots,
\label{eqn-A06}
\end{equation}
where $y_j \in \R$ ($j=0,1,\dots,\lfloor (P-1)/2
\rfloor$) and $0 \le \alpha_j < 1$ ($j=1,2,\dots,\lfloor (P-1)/2
\rfloor$). 

In what follows, we assume $\xi_{k_0} < 0$ for some nonnegative integer
$k_0$ and then prove the following.

\begin{narrow}
{\it Claim: There exists some $b > 0$ such that $\xi_k < - b$ for
  infinitely many $k$'s.}
\end{narrow}

\noindent
If this is true, (\ref{add-eqn-A04}) implies that $x_k < 0$ for a
sufficiently large $k$, which contradicts the fact that $x_k \ge 0$
for all $k = 0,1,\dots$. As a result, for all $k = 0,1,\dots$, $\xi_k$
must be nonnegative, i.e., the statement (c) is true.

We split $\calA \triangleq \{\alpha_j;j=1,2,\dots,\lfloor
(P-1)/2\rfloor\}$ into rational numbers and irrational numbers.  We
then define $\calA_0$ as the set of the rational numbers of
$\calA$. Next we choose an irrational number $\alpha_{j_1}$ from
$\calA \backslash \calA_0$ (if any) and let $\calA_1 = \{\alpha_j \in
\calA \backslash \calA_0; \alpha_j / \alpha_{j_1}~\mbox{is
  rational}\}$. Further we choose an irrational number $\alpha_{j_2}$
from $\calA \backslash (\calA_0 \cup \calA_1)$ (if any) and let
$\calA_2 = \{\alpha_j \in \calA \backslash (\calA_0 \cup \calA_1);
\alpha_j / \alpha_{j_2}~\mbox{is rational}\}$.  Repeating this
procedure, we can obtain $\tilde{P}$ sets, $\calA_j$'s
($j=1,2,\dots,\tilde{P}$), where $\tilde{P}$ may be equal to zero,
i.e., all members of $\calA$ may be rational. Let $\tilde{\alpha}_j$
($j=0,1,\dots,\tilde{P}$) denote some number such that all members of
$\calA_j$ are multiples of $\tilde{\alpha}_j$. Then
$\tilde{\alpha}_j$'s ($j=0,1,\dots,\tilde{P}$) are linearly
independent over the rationals (see Definition~\ref{def-D1}). Note
here that for $n=1,2,\dots$,
\[
\cos (n t)
= T_n(\cos t), \qquad t \in \R,
\]
where $T_n(t)$'s $(n=1,2,\dots)$ denote the Chebyshev polynomials of
the first kind. It thus follows from (\ref{eqn-A06}) that there exist
some polynomial functions $\psi^{(\calA_j)}$'s
$(j=0,1,\dots,\tilde{P})$ on $\R$ such that
\begin{equation}
\xi_k = y_0 +
\psi^{(\calA_0)} \circ \cos(2\pi k \tilde{\alpha}_0) 
+
\sum_{j=1}^{\tilde{P}}
\psi^{(\calA_j)} \circ \cos(2\pi k \tilde{\alpha}_j),
\qquad k=0,1,\dots,
\label{eqn-A11}
\end{equation}
where $\psi^{(\calA_j)} \circ \cos(\cdot)$ denotes a composite
function $\psi^{(\calA_j)}(\cos(\cdot))$ of functions
$\psi^{(\calA_j)}(\cdot)$ and $\cos(\cdot)$. Since $\tilde{\alpha}_0$
is rational, there exists some $g \in \N$ such that
\begin{equation}
\psi^{(\calA_0)} \circ (2\pi (ng+k) \tilde{\alpha}_0) 
=
\psi^{(\calA_0)} \circ \cos(2\pi k \tilde{\alpha}_0),
\qquad \mbox{for all}~k,n =0,1,\dots.
\label{eqn-A12}
\end{equation}
Therefore in the case of $\tilde{P} = 0$, it follows from
(\ref{eqn-A11}) and (\ref{eqn-A12}) that
\[
\xi_{ng+k_0} 
= y_0 +
\psi^{(\calA_0)} \circ \cos(2\pi k_0 \tilde{\alpha}_0) 
=\xi_{k_0} < 0,
\qquad \mbox{for all}~n =0,1,\dots,
\]
which implies the above claim.

We next consider the case of $\tilde{P} \geq 1$. Since
$g\tilde{\alpha}_1$, $g\tilde{\alpha}_2$, \dots,
$g\tilde{\alpha}_{\tilde{P}}$ are linearly independent over the
rationals, it follows from Proposition~\ref{prop-Kro} that for any
$\varepsilon > 0$ and any $\vc{t} \triangleq
(t_1,t_2,\dots,t_{\tilde{P}}) \in \R^{\tilde{P}}$, there exist
integers $n_{\ast}:= n_{\ast}(\varepsilon,\vc{t})$ and
$l_j:=l_j(\varepsilon,\vc{t})$ ($j=1,2,\dots,\tilde{P}$) such that
\[
\left| (n_{\ast} g +k_0) \tilde{\alpha}_j - l_j - t_j \right|
< {\varepsilon \over 2\pi},
\qquad j=1,2,\dots,\tilde{P}.
\]
Thus since $\psi^{(\calA_j)}\circ \cos(2 \pi x)$ is a continuous
function of $x$, there exists some $\delta:=\delta(\varepsilon) > 0$
such that $\lim_{\varepsilon \downarrow 0} \delta = 0$
and
\begin{equation}
| \psi^{(\calA_j)} 
\circ \cos(2\pi (n_{\ast} g +k_0) \tilde{\alpha}_j) 
- 
\psi^{(\calA_j)} 
\circ \cos (2 \pi t_j)  |
<  \delta,
\qquad j=1,2,\dots,\tilde{P}.
\label{eqn-A13}
\end{equation}
It follows from (\ref{eqn-A11}), (\ref{eqn-A12}) and (\ref{eqn-A13})
that
\begin{eqnarray}
\lefteqn{
\left| \xi_{n_{\ast}g + k_0}  - \left( y_0 + 
\psi^{(\calA_0)} \circ \cos (2\pi  k_0 \tilde{\alpha}_0) 
+
\sum_{j=1}^{\tilde{P}} 
\psi^{(\calA_j)} \circ \cos(2  \pi t_j) \right) \right|
}
\quad &&
\nonumber
\\
&& \le \sum_{j=1}^{\tilde{P}}
\left| 
\psi^{(\calA_j)} 
\circ \cos(2\pi (n_{\ast} g +k_0) \tilde{\alpha}_j) 
- 
\psi^{(\calA_j)} \circ \cos(2  \pi t_j)
\right|
< \tilde{P}\delta. \quad
\label{eqn-A14}
\end{eqnarray}
We define $V_{+}(k)$ and $V_{-}(k)$ ($k=0,1,\dots$) as
\begin{eqnarray*}
V_{+}(k) 
&=&  y_0 
+  \psi^{(\calA_0)} \circ \cos(2\pi k \tilde{\alpha}_0)
+
\max_{\svc{t} \in \R^{\tilde{P}}} 
\sum_{j=1}^{\tilde{P}}
\psi^{(\calA_j)}\circ \cos (2 \pi t_j),
\\
V_{-}(k) 
&=&  y_0 
+ \psi^{(\calA_0)} \circ \cos(2\pi k \tilde{\alpha}_0) 
+
\min_{\svc{t} \in \R^{\tilde{P}}} 
\sum_{j=1}^{\tilde{P}} 
\psi^{(\calA_j)} \circ \cos(2  \pi t_j),
\end{eqnarray*}
respectively.  It follows from the above definition and
(\ref{eqn-A11}) that $V_-(k_0) \le \xi_{k_0} \le V_+(k_0)$. Further
(\ref{eqn-A14}) implies that $\{\xi_{ng + k_0};n=0,1,\dots\}$ is
dense in the interval $[V_{-}(k_0), V_{+}(k_0)]$. Thus there exist
infinitely many $n$'s such that $\xi_{ng + k_0} < V_-(k_0)/2 <0$. This
completes the proof of the statement (c).

\medskip
\noindent
\quad {\it Statement (d).} We prove this by reduction to absurdity,
assuming $\xi_{\hat{k}} \le 0$ for some nonnegative integer $\hat{k}$.
Since $(\arg \sigma_j) / \pi $ is a rational number for any
$j=0,1,\dots,P-1$, there exist a positive integer $g$ and nonnegative
integers $l_0,l_1,\dots,l_{P-1}$ such $\sigma_j = \sigma \exp(\iota
2\pi l_j/g)$ ($j=0,1,\dots,P-1$).  Clearly, $\xi_{ng+\hat{k}} \le 0$
for all $n=0,1,\dots$. It thus follows from (\ref{add-eqn-A04}) that
for any $\varepsilon > 0$ there exists some nonnegative integer
$\hat{n}$ such that $\hat{n}g+\hat{k} \ge \breve{m}-1$ and
\[
x_{ng+\hat{k}} 
\le \varepsilon {(ng+\hat{k})^{\breve{m}-1} \over \sigma^{ng+\hat{k}}},
\qquad \mbox{for all}~n=\hat{n},\hat{n}+1,\dots.
\]

We now fix $\hat{n}$ to be such that $\{x_k\}$ is nonincreasing for
all $k \ge \hat{n}g+\hat{k}$ (recall that $\{x_k\}$ is eventually
nonincreasing)\footnote{This part is the only difference from the
  proof of the published version.}. It then follows that
\[
x_{ng+\hat{k}+l} 
\le \varepsilon {(ng+\hat{k})^{\breve{m}-1} \over \sigma^{ng+\hat{k}}},
\qquad n=\hat{n},\hat{n}+1,\dots,~~l = 0,1,\dots,g-1,
\]
which yields for $0 \le y < \sigma$,
\begin{eqnarray}
0 \le f(y)
&\le& \sum_{k=0}^{\hat{n}g + \hat{k}-1} y^k x_k
+ \varepsilon \sum_{n=\hat{n}}^{\infty} 
{(ng+\hat{k})^{\breve{m}-1} \over \sigma^{ng+\hat{k}}} y^{ng+\hat{k}}
\sum_{l=0}^{g-1} y^l
\nonumber
\\
&\le& C_1
+ \varepsilon {1 - \sigma^g \over 1 - \sigma} 
\sum_{n=\hat{n}}^{\infty} (ng+\hat{k})^{\breve{m}-1}
\left({y \over \sigma}\right)^{ng+\hat{k}}
\nonumber
\\
&\le& 
C_1 + \varepsilon C_2 \sum_{k=\breve{m}-1}^{\infty} k^{\breve{m}-1}
\left({y \over \sigma}\right)^k
\nonumber
\\
&\le&  
C_1 + \varepsilon C_2 
\sum_{k=\breve{m}-1}^{\infty} k^{\breve{m}-1}
\left({y \over \sigma}\right)^{k-\breve{m}+1},
\label{eqn-A08}
\end{eqnarray}
where $C_1 = \sum_{k=0}^{\hat{n}g + \hat{k}-1} \sigma^k x_k < \infty$
and $C_2 = (1 - \sigma^g)/(1 - \sigma)$. Note here that the second
last inequality in (\ref{eqn-A08}) follows from $\hat{n}g+\hat{k} \ge
\breve{m}-1$ and the last one follows from $0\le y/\sigma < 1$.

Let $\phi(y) =
\sum_{k=0}^{\infty} (y / \sigma)^k=-\sigma(y-\sigma)^{-1}$ for $0 \le
y < \sigma$. We then have for $0 \le y < \sigma$,
\[
{\rmd^{\breve{m}-1} \over \rmd y^{\breve{m}-1}} \phi(y)
=
\sum_{k=0}^{\infty} {\rmd^{\breve{m}-1} \over \rmd y^{\breve{m}-1}} 
\left({y \over \sigma}\right)^k
= {1 \over \sigma^{\breve{m}-1}}
\sum_{k=\breve{m}-1}^{\infty} k(k-1) \cdots (k - \breve{m} + 2)
\left({y \over \sigma}\right)^{k-\breve{m}+1}.
\]
Thus for $1 \le l \le \breve{m}-1$,
\[
\sum_{k=\breve{m}-1}^{\infty} k(k-1) \cdots (k-l+1)
\left({y \over \sigma}\right)^{k-\breve{m}+1}
\le \sigma^{\breve{m}-1}
{\rmd^{\breve{m}-1} \over \rmd y^{\breve{m}-1}} \phi(y).
\]
Using this inequality and (\ref{add-eqn-A07}), we can bound $f(y)$ in
(\ref{eqn-A08}) as follows.
\[
0 \le
f(y) \le C_1 
+ \varepsilon C {\rmd^{\breve{m}-1} \over \rmd y^{\breve{m}-1}} \phi(y),
\]
where $C = C_2\sigma^{\breve{m}-1}\sum_{l=1}^{\breve{m}-1}b_l$. Further,
\[
{\rmd^{\breve{m}-1} \over \rmd y^{\breve{m}-1}} \phi(y)
= -\sigma {\rmd^{\breve{m}-1} \over \rmd y^{\breve{m}-1}} (y-\sigma)^{-1}
=  \sigma(-1)^{\breve{m}} (\breve{m}-1)! (y-\sigma)^{-\breve{m}}.
\]
As a result,
\[
0\le 
\limsup_{y \uparrow \sigma} 
\left(1 - {y \over \sigma} \right)^{\breve{m}} f(y) 
\le \varepsilon C (\breve{m}-1)! \sigma^{-\breve{m}+1}.
\]
Letting $\varepsilon \to 0$ in the above inequality, we have $\lim_{y
  \uparrow \sigma} \left(1 - {y /\sigma} \right)^{\breve{m}} f(y) =
0$, which contradicts Assumption~\ref{assu-append-1}.  \qed

\subsection{Proof of Lemma \ref{add-lem-period} }\label{proof-lem-add-lem-period} 

Under Assumption \ref{assu-MAdP}, for any $i,j \in \J$ ($i \neq j$)
there exist integers $k_{i,j}$ and $k_{j,i}$ ($k_{i,j}+k_{j,i} \neq
0$) such that $(0,i) \rightarrow (k_{i,j},j)$ and $(0,j) \rightarrow
(k_{j,i},i)$. Let $\K_{j\rightarrow i \rightarrow j}$ denote
\begin{equation}
\K_{j\rightarrow i \rightarrow j}
=
\{
k_{j,i}+k_{i,j}
\}
\cup
\{
k_{j,i} + k + k_{i,j}; \; k \in \K_{i}
\}.
\label{def-K_{jij}}
\end{equation}
Clearly $\K_{j\rightarrow i \rightarrow j} \subseteq \K_j$ and therefore
\begin{equation}
\gcd\{k \in \K_{j\rightarrow i \rightarrow j}\} \ge
\gcd\{k \in \K_j\}= d_j.
\label{eqn-gcd_K_{jij}}
\end{equation}
In what follows, we prove $\gcd\{k \in \K_{j\rightarrow i \rightarrow
  j}\} \le d_i$, from which and (\ref{eqn-gcd_K_{jij}}) it follows
that $d_j \le d_i$. Interchanging $i$ and $j$ in the proof of $d_j \le
d_i$, we can readily show that $d_{i} \le d_{j}$. Therefore we have
$d_{i} = d_{j}$.

Since $(0,i) \rightarrow (k_{i,j},j) \rightarrow (k_{i,j}+k_{j,i},i)$,
we have $k_{i,j}+k_{j,i} \in \K_{i}$ and therefore $k_{i,j}+k_{j,i} =
a_0d_i$ for some integer $a_0 \neq 0$. Note here that $\K_i$ has at
least two elements because
\[
\{
k_{i,j}+k_{j,i}
\}
\cup
\{
k_{i,j} + k + k_{j,i}; \; k \in \K_j
\}
\subseteq \K_i.
\]
Thus there exists a couple of nonzero integers $(a_1,a_2)$ such
that $\{a_1d_i, a_2d_i\} \subseteq \K_i$ and $\gcd\{a_1,a_2\} =~1$,
due to $d_i = \gcd \{k \in \K_i\}$. It follows from
(\ref{def-K_{jij}}) and $k_{i,j}+k_{j,i} = a_0d_i$ that
\begin{eqnarray*}
\K_{j\rightarrow i \rightarrow j}
&\supseteq&
\{
k_{j,i}+k_{i,j}\} \cup \{k_{j,i}+a_1d_i+k_{i,j},~k_{j,i}+a_2d_i+k_{i,j}
\} 
\nonumber
\\
&=& \{a_0d_i\} \cup \{a_0d_i+a_1d_i,~a_0d_i+a_2d_i\} 
\nonumber
\\
&=&
\{a_0d_i, (a_0+a_1)d_i, (a_0+a_2)d_i\},
\end{eqnarray*}
which leads to $\gcd\{k \in
\K_{j\rightarrow i \rightarrow j}\} \le \gcd\{a_0d_i, (a_0+a_1)d_i,
(a_0+a_2)d_i\} =d_i$.  \qed

\subsection{Proof of Theorem~\ref{append-theorem}}\label{proof-append-theorem}

Since the if-part follows from Lemma~\ref{lem-delta-01}, we prove the
only-if part. Let $\vc{V}(\omega)$ denote a $J \times J$ matrix such
that
\[
\vc{V}(\omega) = \diag(\vc{g}(y))^{-1} \vc{\Gamma}^{\ast}(y\omega)
\diag(\vc{g}(y)),
\qquad |\omega| = 1,
\]
where $\diag(\vc{x})$ denotes a diagonal matrix whose $j$th diagonal
element is equal to $[\vc{x}]_j$ for a vector $\vc{x}$. It is easy to
see that $\vc{V}(1)$ is irreducible and stochastic and
$\delta(\vc{V}(\omega)) = \delta(\vc{\Gamma}^{\ast}(y\omega)) =
1$. Let $\vc{f} = (f_j; j\in \J)$ denote a right eigenvector of
$\vc{V}(\omega)$ corresponding to $\delta(\vc{V}(\omega)) = 1$. We
then have for any $n \in \N$, $(\vc{V}(\omega))^n \vc{f} = \vc{f}$ and
thus
\begin{equation}
f_i = \sum_{j \in \J} [(\vc{V}(\omega))^n]_{i,j}f_j
= \sum_{j \in \J} \sum_{k \in \Z} y^k  [\vc{\Gamma}^{(n)}(k)]_{i,j} 
{[\vc{g}(y)]_j \over [\vc{g}(y)]_i} \cdot \omega^k f_j,
\qquad i,j \in \J,
\label{eqn-add-B00}
\end{equation}
where $\{\vc{\Gamma}^{(n)}(k);k\in\Z\}$ is the $n$th-fold convolution
of $\{\vc{\Gamma}(k);k\in\Z\}$ with itself. Note that $\sum_{j \in \J}
\sum_{k \in \Z} y^k[\vc{\Gamma}^{(n)}(k)]_{i,j} [\vc{g}(y)]_j /
    [\vc{g}(y)]_i=1$ because $(\vc{V}(1))^n\vc{e} = \vc{e}$. Let $i'$
    denote an element of $\J$ such that $|f_{i'}| \ge |f_j|$ for all
    $j \in \J$. It then follows from (\ref{eqn-add-B00}) that for any
    $j \in \J$,
\begin{equation}
\omega^k {f_j \over f_{i'}} = 1~\mbox{if}~[\vc{\Gamma}^{(n)}(k)]_{i',j} > 0.
\label{eqn-add-B01}
\end{equation}
Since $\vc{\Gamma}^{\ast}(1)$ is irreducible, for any
$j \in \J$ there exist some $n \in \N$ and $k \in \Z$ such that
$[\vc{\Gamma}^{(n)}(k)]_{i',j} > 0$. Thus (\ref{eqn-add-B01}) implies
that $|f_j|$'s are all equal, because $|\omega| = 1$. We now consider a
path from phase $i$ to phase $i$ such that
\[
(0,i) \rightarrow (k_1,i_1) \rightarrow (k_2,i_2)
\rightarrow (k_m,i_m) \triangleq (k_m,i),
\]
where $(k_l,i_l) \in \Z \times \J$ for $l=1,2,\dots,m$ and $m \in
\N$. Since the period of MAdP $\{\vc{\Gamma}(k)\}$ is equal to $d$,
$k_1+k_2+\cdots+k_m$ is a multiple of $d$. From (\ref{eqn-add-B01}),
we have $\omega^{k_1+k_2+\cdots+k_m} = 1$ and thus $\omega^d = 1$. The
proof of the only-if part is completed. 

As for the remaining statements, (\ref{add-eqn-52}) is obvious, and it
follows from Lemma~\ref{lem-delta-01}~(a) that if
$\delta(\vc{\Gamma}^{\ast}(y\omega)) = 1$, then the eigenvalue
$\delta(\vc{\Gamma}^{\ast}(y\omega_d^{\nu})) = 1$ is simple for
$\nu=0,1,\dots,d-1$. \qed

\section{Kronecker's Approximation Theorem}\label{appendix-lemmas}

The following is Kronecker's approximation theorem. For details, see,
e.g., Theorem~7.10 in \cite{Apos97}.
\begin{prop}
\label{prop-Kro}
Let $\gamma_i$'s ($i=1,2,\dots,n$) denote arbitrary real numbers. Let
$\beta_i$'s ($i=1,2,\dots,n$) denote arbitrary real numbers such that
$\beta_1,\beta_2,\dots,\beta_n$ and 1 are linearly independent over
the rationals (see Definition~\ref{def-D1} below). Then for any
$\varepsilon > 0$, there exist an $(n+1)$-tuple
$(k,l_1,l_2,\dots,l_n)$ of integers such that
\begin{equation}
\left| k \beta_i - l_i - \gamma_i \right| < \varepsilon, 
\qquad \mbox{for all}~i=1,2,\dots,n,
\label{eq-D02}
\end{equation}
and thus for any $\varepsilon > 0$ and any $\tilde{\gamma}_i \in [0,1]$,
\[
\left| k \beta_i - \lfloor k\beta_i \rfloor - \tilde{\gamma}_i \right| < \varepsilon, 
\qquad \mbox{for all}~i=1,2,\dots,n,
\]
which implies that $k \beta_i - \lfloor k\beta_i \rfloor$ 
$(k \in \Z)$ is dense in the interval $[0,1]$. 
\end{prop}

\begin{defn}\label{def-D1}
Arbitrary real numbers $\beta_i$'s $(i=1,2,\dots,n)$ are said to be
{\it linearly independent over the rationals (equivalently integers)}
if there exists no set of rational numbers $q_i$'s
$(i=1,2,\dots,n)$ such that $(q_1,q_2,\dots,q_n) \neq \vc{0}$ and
\begin{equation}
\beta_1 q_1 + \beta_2 q_2 + \cdots + \beta_n q_n = 0.
\label{add-eqn-D01}
\end{equation}
Therefore if $\beta_i$'s are linearly independent over the rationals,
(\ref{add-eqn-D01}) implies that $q_1 = q_2 = \cdots = q_n = 0$.
\end{defn}

\section{Example against Assumption \ref{assu-3}}\label{appendix-theta}

We suppose $\vc{A}(k)$'s ($k=0,1,\dots$) are scalars such that for
some finite $r > 1$,
\[
\vc{A}(k) 
= \left.{1 \over r^k}{1 \over (k+1)^3} 
\right/ \sum_{n=0}^{\infty}{1 \over r^n}{1 \over (n+1)^3},
\qquad k=0,1,\dots.
\]
Clearly, $r_A = r$. We define $F(x)$ ($x\ge1$) as
\[
F(x) = x\sum_{k=0}^{\infty}{1 \over x^k}{1 \over (k+1)^3}.
\]
It then follows that for any $x \ge 1$,
\[
F'(x) = 1 - \sum_{k=2}^{\infty}{1 \over x^k}{k-1 \over (k+1)^3},
\]
which leads to
\begin{eqnarray*}
F'(1)
&\ge& 1 - \sum_{k=2}^{\infty}{1 \over (k+1)^2} 
= {9 \over 4} - {\pi^2 \over 6} > 0,
\\
F''(x) &>& 0.
\end{eqnarray*}
Thus since $F(r_A) > F(1)$, we have
\[
{1 \over r_A}\delta(\vc{A}^{\ast}(r_A)) = 
{1 \over r_A}\vc{A}^{\ast}(r_A) = {F(1) \over F(r_A)} < 1.
\]

\section*{Acknowledgments}
The authors are thankful to Takine Tetsuya and Masakiyo Miyazawa for
his valuable comments on an early version of this paper. They 
also thank the anonymous referees and the editor for their
suggestions on how to strengthen the results and improve the
presentation of this paper. Research of the third author was supported
in part by Grant-in-Aid for Young Scientists (B) of Japan Society for
the Promotion of Science under Grant No.~21710151.

\renewcommand{\baselinestretch}{0.9}

\end{document}